\documentclass[11pt]{amsart}
\usepackage{amsmath,amsfonts,amscd,amssymb}

\makeatletter \@addtoreset{section}{part} \@addtoreset{subsection}{section}
\makeatother

\newtheorem{theorem}{Theorem}[section]
\newtheorem{proposition}[theorem]{Proposition}
\newtheorem{lemma}[theorem]{Lemma}
\newtheorem{corollary}[theorem]{Corollary}

\theoremstyle{definition}

\theoremstyle{remark}\newtheorem{remark}[theorem]{Remark}

\numberwithin{equation}{section}

\newcommand{\field}[1]{\ensuremath{\mathbb{#1}}}

\newcommand{\CC}{\field{C}}
\newcommand{\HH}{\field{H}}
\newcommand{\NN}{\field{N}}

\newcommand{\PP}{\field{P}}
\newcommand{\RR}{\field{R}}

\newcommand{\ZZ}{\field{Z}}
\newcommand{\J}{\mathbf{J}}

\newcommand{\btau}{\boldsymbol{\tau}}

 \DeclareMathOperator{\Hom}{Hom}

\DeclareMathOperator{\im}{Im} \DeclareMathOperator{\Alb}{Alb}
\DeclareMathOperator{\Bires}{Bires}
\DeclareMathOperator{\Pic}{Pic}
\DeclareMathOperator{\re}{Re}\DeclareMathOperator{\Div}{Div}
\DeclareMathOperator{\Diff}{Diff}\DeclareMathOperator{\Imm}{Im}
\DeclareMathOperator{\Homeo}{Homeo}\DeclareMathOperator{\Mob}{M\ddot{o}b}
\DeclareMathOperator{\Area}{Area}\DeclareMathOperator{\PDiv}{PDiv}
\begin{document}

\title[Free Bosons and Tau-Functions]
{Free Bosons and Tau-Functions for Compact Riemann Surfaces and Closed Smooth Jordan
Curves I. Current Correlation Functions}
\author{Leon A. Takhtajan}
\address{Department of Mathematics \& Department of Physics, Yale University, New Haven
CT 06520-8283 \& 06511-8167
 USA\footnote{On leave of absence from Department of Mathematics SUNY at
Stony Brook, Stony Brook NY 11794-3651 USA.}.}\email{leon.takhtajan@yale.edu \&
leontak@math.sunysb.edu}

\begin{abstract}
We study families of quantum field theories of free bosons on a compact Riemann
surface of genus $g$. For the case $g>0$ these theories are parameterized by
holomorphic line bundles of degree $g-1$, and for the case $g=0$ --- by smooth closed
Jordan curves on the complex plane. In both cases we define a notion of
$\tau$-function as a partition function of the theory and evaluate it explicitly. For
the case $g>0$ the $\tau$-function is an analytic torsion~\cite{AG-B-M-N-V}, and for
the case $g=0$ --- the regularized energy of a certain natural pseudo-measure on the
interior domain of a closed curve. For these cases we rigorously prove the Ward
identities for the current correlation functions and determine them explicitly. For
the case $g>0$ these functions coincide with those obtained in~\cite{K-N-T-Y,R2} using
bosonization. For the case $g=0$ the $\tau$-function we have defined coincides with
the $\tau$-function introduced in~\cite{MWZ,W-Z,K-K-MW-W-Z} as a dispersionless limit of the
Sato's $\tau$-function for the two-dimensional Toda hierarchy. As a corollary of the
Ward identities, we obtain recent results~\cite{W-Z,K-K-MW-W-Z} on relations between
conformal maps of exterior domains and $\tau$-functions. For this case we also define
a Hermitian metric on the space of all contours of given area. As another corollary of
the Ward identities we prove that the introduced metric is K\"{a}hler and the
logarithm of the $\tau$-function is its K\"{a}hler potential.
\end{abstract}

\maketitle \tableofcontents
\section*{Introduction}

The general notion of a Riemann surface was defined by Riemann. He proved that, in
modern terminology, every compact Riemann surface is a Riemann surface of an algebraic
function $f(x,y)=0$. Continuing the work of Abel and Jacobi, Riemann also introduced
the notion of a general theta-function and solved the Jacobi inversion problem. After
Riemann's death, his approach was repeatedly criticized by Weierstrass for the lack of
rigor: in the proof of Dirichlet principle Riemann assumed that Dirichlet functional
has a minimum. Later, Riemann's main results were proved by Schwarz and C.~Neumann
without using this assumption and in 1901 Hilbert finally proved the Dirichlet
principle itself. At the same time Weierstrass criticism stimulated the development of
the algebraic theory of algebraic functions, i.e.~a theory that does not use
complex-analytic methods. This theory was formulated by Brill, M.~Noether, Dedekind
and H.~Weber, quite in parallel with the theory of algebraic numbers. In 1924, Artin
was the first to consider algebraic functions over fields different from the field
$\CC$ of complex numbers. Artin's results were generalized by F.K.~Schmidt who, in
particular, developed the theory of algebraic functions over a finite field. This
development finally led A.~Weil, among other things, to the general formulation of the
algebraic geometry over an arbitrary field (see~\cite{iwasawa} for an exposition and
references).

This old story nowadays is again in the spotlight due to the advent of the string
theory. It is amazing that dramatic applications of quantum fields and strings to
various mathematical areas are quite similar in spirit to Riemann's original work.
Thus in his proof of the existence of a harmonic function with prescribed
singularities at given points, Riemann used the analogy with electrostatic theory,
assuming as obvious that for any charge distribution there exists an electrostatic
potential, or equivalently, there exists a flow of the ideal fluid with prescribed
sources and sinks (see, e.g.,\cite{C-H}). Using modern terminology, Riemann was
exploiting methods of classical field theory for mathematical purposes, ``probing''
mathematical objects with physical theories and translating the physical ``output''
back into the mathematical statements. The same exact idea is in the heart of today's
applications, with the ``only'' difference that classical fields are being replaced by
quantum fields and strings. Succinctly, this idea can be described as follows. In
classical field theory one studies critical points of action functionals (like
Dirichlet functional) which satisfy partial differential equations. In quantum field
theory one studies partition and correlation functions of quantum fields defined by
Feynman path integrals. In this formalism the critical values of the action functional
as well as ``higher invariants'' --- so-called quantum corrections, naturally appear
in the perturbative expansion (see~\cite{QFT} for exposition aimed at mathematicians).
When probing mathematical objects (e.g.~topological, smooth or complex-analytic
manifolds) by quantum field theories, the mathematical output is encoded in partition
and correlation functions that are expanded in terms of the critical values and
quantum corrections.

Symmetries of physical theories play a fundamental role. Continuous symmetries in
classical field theories correspond to conservation laws, and in quantum field
theories they correspond to the so-called Ward identities expressing correlation
functions through the partition function. The exploitation of symmetry is a powerful
tool for studying conformal field theories in two dimensions, as discovered by
A.B.~Belavin, A.M.~Polyakov and A.B.~Zamolodchikov~\cite{BPZ}. In addition to
important physical applications to the critical phenomena in statistical mechanics,
conformal field theories have a rich mathematical structure with applications to
representation theory, complex analysis, uniformization of Riemann surfaces and
complex algebraic geometry of moduli spaces.

Quantum field theories of free fermions and free bosons on compact Riemann surfaces
are basic examples of conformal field theories. It was proved by L.~Alvarez-Gaum\'{e},
J.-B.~Bost, G.~Moore, P.~Nelson, and C.~Vafa and by E.~Verlinde and H.~Verlinde
(see~\cite{AG-M-V,AG-B-M-N-V,B-N,V-V} and references therein) that these theories are
equivalent. Namely, there exists a remarkable correspondence between bosons and
fermions, called bosonization, that expresses partition and correlation functions for
one theory in terms of the other theory. Since for both theories these functions can
be written down explicitly in terms of theta-functions and related algebro-geometric
objects, bosonization yields non-trivial identities between them. In particular, it
gives another proof of the celebrated Fay's trisecant identity for
theta-functions~\cite{Fay0}, which play a fundamental role in the KP hierarchy (see,
e.g.~\cite{Mumford}). Approach in~\cite{AG-B-M-N-V,V-V} was based on path integrals
with mathematical proofs using analytic torsion and Quillen's type isometries of
determinant line bundles for $\bar\partial$-operators. Another approach to
bosonization based on Sato's infinite Grassmannian manifold was developed by
L.~Alvarez-Gaum\'{e}, C.~Gomez and C.~Reina~\cite{AG-G-R}, C.Vafa~\cite{V} and by
N.~Kawamoto, Y.~Namikawa, A.~Tsuchiya and Y.~Yamada~\cite{K-N-T-Y}. Purely
algebro-geometric approach was given by A.~Raina~\cite{R1,R2}. In addition to the
proof of Fay's trisecant identity, in papers~\cite{K-N-T-Y,R2} bosonization was also
used to compute explicitly bosonic current correlation functions from correlation
functions for fermion operators.

This story brings us to Part 1 of the paper, where we derive bosonic current
correlation functions directly from the basic $U(1)$-symmetry of the theory. We define
$\btau$-function as a partition function and, following ~\cite{AG-B-M-N-V}, show that
it is essentially an absolute value of Sato's $\tau$-function computed
in~\cite{K-N-T-Y}. Using explicit description of the complex structure of the Jacobian
variety of a compact Riemann surface and basic properties of the Abel-Jacobi map we
prove $U(1)$-gauge symmetry Ward identities. As the result, we completely determine
all reduced normalized multi-point current correlation functions in terms of partial
derivatives of $\log\btau$ with respect to the complex coordinates on the Jacobian.
Obtained formulas are in perfect agreement with results in~\cite{K-N-T-Y,R2}. Though
our exposition is based on path integrals, all integrals that actually appear are
Gaussian and have pure mathematical definition. Therefore our approach is completely
rigorous.

It is well-known that Ward identities for other quantum field theories also encode
important mathematical information. An interesting example is given by the quantum
Liouville theory --- the theory that arises as a conformal anomaly in Polyakov's
approach to the string theory~\cite{Pol2}. It was conjectured by Polyakov~\cite{P3}
that the semi-classical limit of conformal Ward identities for the Liouville theory
yields an expression for the Poincar\'{e}'s accessory parameters in terms of the
critical value of  Liouville action functional. Together with P.G.~Zograf we have
proved this and other results on the relation between accessory parameters and the
Weil-Petersson geometry of the Teichm\'{u}ller space~\cite{Z-T1,Z-T2,Z-T3}.
In~\cite{Tak} we summarized the geometric approach to two-dimensional quantum gravity
and interpreted results in~\cite{Z-T1,Z-T2,Z-T3} as conformal Ward identities for
multi-point correlation functions of stress-energy tensor components.

From this perspective, Part 1 of this paper is just another application of the same
idea. The only difference is that instead of the quantum Liouville theory we consider
simpler quantum theory of free bosons on a Riemann surface and instead of correlation
functions with stress-energy tensor components we consider current correlation
functions, which are also much simpler. Also, instead of the ``nonlinear''
Teichm\"{u}ller theory --- the deformation theory of complex structures on a given
compact Riemann surface used in~\cite{Z-T1,Z-T2,Z-T3}, we are using a ``linear
theory'' of the Jacobian variety --- the deformation theory of degree zero holomorphic
line bundles over the Riemann surface.

In Part 2, we exploit the same simple idea for the case of free bosons on bounded
simply-connected domains in the complex plane $\CC$. Using results of Part 1 as a
motivation, for every smooth closed Jordan curve $C$ in $\CC$ we formulate a quantum
field theory on the Riemann sphere $\PP^1$ that depends on the interior domain
$\Omega$ of a contour $C$. We define the $\tau$-function of a contour $C$ as a
normalized partition function of the theory and compute it explicitly. It turns out to
be the exponential of the regularized energy of a certain pseudo-measure on the domain
$\Omega$ and, quite remarkably, coincides with the Mineev-Weinstein--Wiegmann--Zabrodin
$\tau$-function~\cite{W-Z,K-K-MW-W-Z}, defined as a dispersionless limit of the Sato's
$\tau$-function for the two-dimensional Toda hierarchy. We prove the Ward identities
which express through $\log\tau$ the difference between current correlation functions
of free bosons parameterized by $C$ and current correlation functions of free bosons
on the exterior domain $\PP^1\setminus\Omega$ satisfying the Dirichlet boundary
condition. Correlation functions for free bosons on $\PP^1\setminus\Omega$ with the
Dirichlet boundary condition can be written explicitly in terms of the conformal map
$G$ of the domain $\PP^1\setminus\Omega$ onto the exterior of the unit disk. Thus the
Ward identities give another proof of remarkable relations between the conformal map
$G$ and the $\tau$-function, discovered by M.~Mineev-Weinstein, P.B.~Wiegmann and
A.~Zabrodin~\cite{MWZ} in the content of the physical problem of 2D interface
dynamics and studied extensively in~\cite{W-Z,K-K-MW-W-Z}. We
also introduce a Hermitian metric $H$ on the space of all contours of given area and
prove as a corollary of the Ward identities that this metric is K\"{a}hler and the
logarithm of the $\tau$-function is a K\"{a}hler potential of the metric $H$.

The existence of the conformal map $G$ is a consequence of the  Riemann mapping
theorem that every simply-connected domain in the complex plane $\CC$ whose boundary
consists of more then two points can be conformally mapped onto the unit disk. There
is a Riemann's proof of this theorem which is based on the classical field theory and
defines the conformal mapping through the complex potential of a certain charge
distribution/fluid flow in the domain~\cite{C-H}. Corresponding quantum field theory
naturally introduces the $\tau$-function of contour $C$. Just as conformal Ward
identities for the quantum Liouville theory imply that the critical value of the
Liouville action is a generating function for the accessory
parameters~\cite{Z-T1,Z-T2,Z-T3}, current Ward identities for the theory of free
bosons on $\PP^1$ imply that the $\tau$-function of smooth bounded contours is a
generating function for the conformal maps. The role of Teichm\"{u}ller theory in this
case is played by the deformation theory of contours as it has been stated by
I.M.~Krichever~\cite{K}.

Thus both parts of the paper are related not just by exploiting the same idea of Ward
identities for current correlation functions, but also by using Riemann's two major
achievements in complex analysis as a framework for quantum field theories.

The above examples illustrate the importance of probing one-dimensional complex
manifolds: compact Riemann surfaces and simply-connected domains in the complex plane,
by quantum field theories. This idea can be pushed further by considering adelic
formulation of quantum field theories on algebraic curves over an arbitrary field of
constants~\cite{W1,W2,T} where Ward identities result in reciprocity laws. There is
also a possibility of defining quantum field theories on the fields of algebraic
numbers and it seems that the proper setting should be ``quantum field theories on
regular one-dimensional schemes''.

In this paper we are systematically using physical terminology. This is done for the
only purpose: to show the true origin of the ideas and methods, and not for the
purpose of proofs. Mathematically oriented reader can completely ignore these terms,
which are \emph{emphasized} in the main text on their first appearance, and consider
them simply words\footnote{``Guaranteed to raise a smile'' from physicists.}. All
results in the paper are rigorously proved (with usual space limitations) and no
knowledge of quantum physics is required for what follows.

Here is a more detailed description of the content of the paper. Section 1 of Part 1
reviews some necessary mathematical facts: Jacobians, properties of the Abel-Jacobi
map and theta-functions in 1.1, Green's functions 1.2. In addition to standard
definition and properties of the Green's function of a Laplace operator on an compact
Riemann surface $X$ of genus $g$, in Lemma 1.1 we included a proof of the Fay's
formula relating two classical kernels on $X$. In Section 2 we define the theory of
free $U(1)$-bosons on $X$ with classical fields\footnote{It should be always clear
from the content whether $g$ is a genus, $U(1)$-valued bosonic field, or an inverse to
a conformal map $G$ as in Part 2.} $g\in C^\infty (X,U(1))$ and field currents
$\J:=g^{-1}d g/2\pi\sqrt{-1}$, parameterized by the holomorphic line bundle $L$ over
$X$ of degree $g-1$ with $h^0(L)=0$. We introduce the action functional $S(\J)$ of the
theory as a sum of the action functional for the standard free bosons on $X$ and the
topological term depending on $L$. We follow~\cite{AG-B-M-N-V} with minor technical
improvement in the invariant definition of the topological term. We define
$\btau$-function --- the partition function as
\begin{equation*}
\btau=\langle\mathbf{1}\rangle_L:=\int_{C^\infty(X,U(1))/U(1)}
[\mathcal{D}\J]\,e^{-2\pi S_L(\J)}
\end{equation*}
and in Proposition 2.3 evaluate it explicitly in terms of the Riemann theta-function.
The computation is a slight simplification of the original proof~\cite{AG-B-M-N-V}.
Namely, we apply the Poisson summation formula for the lattice $\ZZ^{2g}$ instead of
the lattice $\ZZ^g$ as in~\cite{AG-B-M-N-V} and use the transformation formula for the
Riemann theta-function. In Section 3 we prove $U(1)$-gauged symmetry Ward identities
for normalized reduced current correlation functions. Specifically, we introduce
holomorphic and anti--holomorphic components of the bosonic field current as
$J=-g^{-1}\partial g$ and $\bar J=-g^{-1}\bar\partial g$, and define multi-point
current correlation functions as
\begin{multline*}
\langle J(P_1)\dotsb J(P_m)\bar{J}(Q_1)\dotsb\bar{J}(Q_n)\rangle \\
:=\int_{C^\infty(X,U(1))/U(1)}[\mathcal{D}\J] \,J(P_1)\dotsb
J(P_m)\bar{J}(Q_1)\dotsb\bar{J}(Q_n)\,e^{-2\pi S_L(\J)}.
\end{multline*}
The $U(1)$-gauged Ward identities for normalized 1-point correlation function have the
form
\begin{equation*}
\frac{\partial\log\btau}{\partial z_i}= \int_{a_i}\langle\langle
J(P)\rangle\rangle,\,i=1,\dotsc,g,
\end{equation*}
where $z_i$ are complex coordinates of the point $Z$ on the Jacobian $J(X)$
corresponding to the line bundle $L$, and directly follow from Riemann bilinear
relations.

The Fay's formula is relevant for the computation of the reduced normalized 2-point
correlation function. General expression for the multi-point correlation functions is
given in Theorem 3.3 and is in perfect agreement with~\cite{K-N-T-Y,R2}.

Basically, Part 1 serves as a motivation for Part 2. In Section 1 of Part 2 we recall
necessary mathematical facts. In Section 1.1 we, following
A.A.~Kirillov~\cite{Kirillov1}, review the infinite-dimensional space $\mathcal{C}$ of
all smooth closed Jordan curves on $\CC$ encircling the origin $0$, and other spaces
related to it. We also introduce a double $\PP^1_C$ of the exterior domain
$\PP^1\setminus\Omega$ of the contour $C$. In Section 1.2 we recall basic facts about
classical Green's functions on $\PP^1$ and on $\PP^1\setminus\Omega$, including
Schiffer kernel and Bergman's reproducing kernel. Section 1.3 is devoted to detailed
exposition of Krichever's deformation theory, outlined in~\cite{K}. We introduce the
analogs of Faber polynomials for the conformal map $G$ and define the harmonic moments
of interior $t_0, t_n$ of the contour $C$ and harmonic moments of exterior $v_0, v_n$
of $C$, $n\in\NN$. We formulate and prove Theorem 1.5 which contains explicit
description of holomorphic tangent and cotangent spaces to the infinite-dimensional
complex manifold $\Tilde{\mathcal{C}}_a$ of contours of fixed area $a>0$. The main
ingredient is the so-called Krichever's lemma~\cite{K} that expresses vector fields
$\partial/\partial t_0,
\partial/\partial t_n, \partial/\partial \bar{t}_n$ on $\mathcal{C}$ as meromorphic
$(1,0)$ - forms on $\PP^1_C$. In addition to~\cite{K}, we introduce natural Hermitian
metric $H$ on the spaces $\Tilde{\mathcal{C}}_a$ using the Bergmann reproducing
kernel. The metric $H$ turns out to be K\"{a}hler (see Section 3 and description
below).

 In Section 2 we define theories of free bosons on $\PP^1$ parameterized by
smooth Jordan contours $C\in\mathcal{C}$. The action functional $S_C(\varphi)$ of the
theory is defined as a sum of the action functional for the standard free bosons on
$\PP^1$ and of the analog of the topological term depending on $C$. The partition
function is defined as
\begin{equation*}
\langle\mathbf{1}\rangle_C:=\int_{C^\infty(\PP^1,\RR)/\RR}
[\mathcal{D}\varphi]\,e^{-\frac{1}{\pi}S_C(\varphi)},
\end{equation*}
and the $\tau$-function of the contour $C$ is given by
\begin{equation*}
\tau:=\frac{\langle\mathbf{1}\rangle_C}{\langle\mathbf{1}\rangle_\emptyset},
\end{equation*}
where $\emptyset$ is the empty set that formally corresponds to the case when no
contour is present. In proposition 2.1 we compute the $\tau$-function explicitly. It
turns out to be the exponential of a regularized energy of the pseudo-measure given by
the difference between the Lebesgue measure on the domain $\Omega$ and the
delta-measure at $0$ times the Euclidean area of $\Omega$ and coincides with
Mineev-Weinstein--Wiegmann--Zabrodin $\tau$-function~\cite{MWZ,W-Z,K-K-MW-W-Z}. In Section
3 we introduce
holomorphic and anti-holomorphic components $\jmath=\partial\varphi,
\bar\jmath=\bar{\partial}\varphi$ of the field current $d\varphi$ and define
multi-point current correlation functions
\begin{multline*}
\langle \jmath(z_1)\dotsb \jmath(z_m)\bar{\jmath}(w_1)\dotsb\bar{\jmath}(w_n)\rangle
\\ :=\int_{C^\infty(\PP^1,\RR)/\RR}[\mathcal{D}\varphi] \,\jmath(z_1)\dotsb
\jmath(z_m)\bar{\jmath}(w_1)\dotsb\bar{\jmath}(w_n)\,e^{-\frac{1}{\pi}S_C(\varphi)}
\end{multline*}
for free bosons on $\PP^1$ parameterized by $C\in\mathcal{C}$. In a similar fashion we
introduce current correlation functions for free bosons on $\PP^1\setminus\Omega$ with
the Dirichlet boundary condition, denoted by $\langle\,\,\,\,\,\,\rangle_{DBC}$. In
Theorem 3.1 we prove the Ward identity for the 1-point normalized correlation function
\begin{equation*}
\frac{\partial\log\tau}{\partial t_n}=-\frac{1}{2\pi i} \int_C\langle\langle\jmath(z)
\rangle\rangle z^n dz,\,n\in\NN.
\end{equation*}
As a corollary we immediately get that the $\tau$-function is a generating function
for the harmonic moments of interior, as first proved in~\cite{MWZ,W-Z,K-K-MW-W-Z,K},
\begin{equation*}
v_0=\frac{\partial\log\tau}{\partial t_0},~v_n=\frac{\partial\log\tau}{\partial t_n},
\,n\in\NN.
\end{equation*}
In Corollary 3.7 we derive another Mineev-Weinstein--Wiegmann--Zabrodin result~\cite{MWZ,W-Z}
--- an ``explicit formula'' for the conformal map $G$
\begin{equation*}
\log G(z) =\log z - \frac{1}{2}\frac{\partial^2\log\tau}{\partial t_0^2} -
\sum_{n=1}^{\infty}\frac{z^{-n}}{n}\frac{\partial^2\log\tau}{\partial t_0\partial
t_n}.
\end{equation*}
In theorem 3.9 we prove the following Ward identities for normalized reduced 2-point
correlation functions
\begin{align*}
\frac{\partial^2\log\tau}{\partial t_m\partial t_n}&=\frac{1}{(2\pi i)^2}\int_C\int_C
\bigg(\frac{G^\prime(z)G^\prime(w)}{(G(z)-G(w))^2}-\frac{1}{(z-w)^2}\bigg)z^m w^n dz
dw \\ &=\frac{1}{(2\pi i)^2}\int_C\int_C z^m
w^n\left(\langle\langle\jmath(z)\jmath(w)\rangle\rangle -
\langle\langle\jmath(z)\jmath(w)\rangle\rangle_{DBC} \right)
\end{align*}
and
\begin{align*}
\frac{\partial^2\log\tau}{\partial t_m\partial \bar{t}_n} & = -\frac{1}{(2\pi
i)^2}\int_{C_+}\int_{C_+}
\frac{G^\prime(z)\overline{G^\prime(w)}}{(1-G(z)\overline{G(w)})^2}z^m \bar{w}^n dz
d\bar{w} \\ & = -\frac{1}{(2\pi i)^2}\int_{C_+}\int_{C_+} z^m \bar{w}^n
\left(\langle\langle\jmath(z)\bar{\jmath}(w)\rangle\rangle -
\langle\langle\jmath(z)\bar{\jmath}(w)\rangle\rangle_{DBC} \right),
\end{align*}
where $C_+$ is an arbitrary contour around 0 containing the contour $C$ inside.

The latter formula shows that the natural Hermitian metric $H$ on infinite dimensional
complex manifolds $\Tilde{\mathcal{C}}_a$ is K\"{a}hler and $\log\tau$ is its
K\"{a}hler potential.This is a new result.

In numerous remarks throughout the paper we point out to interesting connections with
other fields, classical complex analysis and theory of the univalent functions in
particular. We plan to return to these questions, as well as to the discussion of the
bosonization, in a separate paper.

\subsection*{Acknowledgments} I very much appreciate stimulating discussions with
I.M.~Krichever, P.B.~Wiegmann and A.~Zabrodin that led to Part 2 of the paper. I would
also like to thank D.~Ben-Zvi for the discussion of current correlations functions on
compact Riemann surfaces. This research was partially supported by the NSF grant
DMS-9802574.

\part{Free Bosons and Tau-Functions for Compact Riemann Surfaces}
\section{Mathematical set-up}
Here we recall, in a succinct form, necessary facts from complex (algebraic) geometry
of compact Riemann surfaces (algebraic curves). The standard reference is~\cite{G-H};
see also~\cite{M,F-K} as well as~\cite{AG-B-M-N-V} for a ``crash course'' for
physicists.

Let $X$ be a compact Riemann surface of genus $g$, which is always assumed to be
connected and without the boundary. Denote by $K=K_X$ the canonical line bundle on $X$
--- holomorphic cotangent bundle of $X$, by $\Div(X)$ --- the group of divisors on $X$,
and by $\Pic(X)$
--- the Picard group of isomorphism classes of holomorphic line bundles over $X$. The
correspondence between line bundles and divisors provides canonical isomorphism
$\Pic(X)\simeq\Div(X)/\PDiv(X)$ --- the group of divisors on $X$ modulo the subgroup
$\PDiv(X)$ of principal divisors. Let $\check{H}^i(X,\mathcal{F})$ be the \v{C}ech
cohomology groups with coefficients in a sheaf $\mathcal{F}$ on $X$ and set
$h^i(\mathcal{F})=\dim_{\CC}\check{H}^i(X,\mathcal{F}),~i=0,1$. As usual, we denote by
$L$ a holomorphic line bundle over $X$ as well as the sheaf of germs of holomorphic
sections of $L$. The Riemann-Roch formula is
\begin{equation*}
h^0(L)-h^0(K\otimes L^{-1})=\deg L + 1-g,
\end{equation*}
where $\deg L$ is a degree of the line bundle $L$.

\subsection{Jacobians and theta-functions} The period map
\begin{equation*}
\check{H}^0(X,K)\ni\omega \mapsto\int_c\omega\in\CC,
\end{equation*}
for all homology classes $[c]\in H_1(X,\ZZ)$ of $1$-cycles $c$ on $X$, defines
canonical inclusion of $H_1(X,\ZZ)$ into $\check{H}^0(X,K)^\vee$ --- the dual vector
space to $\check{H}^0(X,K)$. The additive subgroup $H_1(X,\ZZ)$ of
$\check{H}^0(X,K)^\vee$ is a discrete subgroup of maximal rank $2g$ over $\RR$. The
Albanese variety $\Alb(X)$ of a compact Riemann surface $X$ is canonically defined as
the following $g$-dimensional complex torus
\begin{equation*}
\Alb(X):=\check{H}^0(X,K)^\vee/H_1(X,\ZZ).
\end{equation*}
The Albanese variety $\Alb(X)$ is a complex projective manifold. It carries
translation-invariant K\"{a}hler metric, defined by the following Hermitian inner
product at $\check{H}^0(X,K)=T^{\ast}_0\Alb(X)$ --- the holomorphic cotangent vector
space to $\Alb(X)$ at $0$,
\begin{equation*}
<\omega_1,\omega_2>:=\frac{\sqrt{-1}}{2}\int_X \omega_1\wedge\bar{\omega}_2,~
\omega_1, \omega_2\in\check{H}^0(X,K).
\end{equation*}

Let $\Pic^0(X)$ be the group of degree $0$ line bundles over $X$ --- the identity
component of the Picard group $\Pic(X)$. It follows from the standard exponential
exact sequence and from the Dolbeault isomorphism
\begin{equation*}
\Pic^0(X)=\check{H}^1(X,\mathcal{O})/\check{H}^1(X,\ZZ)=H^{0,1}_{\bar{\partial}}
(X)/H^1(X,\ZZ),
\end{equation*}
where $\mathcal{O}$ is the structure sheaf on $X$ and
$H^1(X,\ZZ):=\Hom(H_1(X,\ZZ),\ZZ)$ is a lattice in the de Rham cohomology group
$H^1(X,\RR)$ of cohomology classes of 1-forms with integral periods. By the Hodge
theorem, $H^1(X,\ZZ)\simeq\mathcal{H}^1(X,\ZZ)$ --- the space of harmonic 1-forms on
$X$ with integral periods, and the mapping
\begin{equation*}
\mathcal{H}^1(X,\ZZ)
\ni\omega=\overline{\omega^{0,1}}+\omega^{0,1}\mapsto\omega^{0,1}\in
H^{0,1}_{\bar{\partial}}(X)
\end{equation*}
gives canonical embedding $H^1(X,\ZZ)\hookrightarrow H^{0,1}_{\bar{\partial}}(X)$. The
following canonical isomorphisms
\begin{equation*}
\Div^0(X)/\PDiv(X)\simeq\Pic^0(X)\simeq\Alb(X)
\end{equation*}
follow from the Serre and Poincar\'{e} dualities.

A Riemann surface $X$ is called Torelli marked if it is equipped with a symplectic
basis $\{a_i,b_i\}_{i=1}^g$ for $H_1(X,\ZZ)$. For a Torelli marked Riemann surface
$X$, let $\{\omega_i\}_{i=1}^g$ be the basis for $\check{H}^0(X,K)$ of abelian
differentials with normalized $a$-periods: $\int_{a_i}\omega_j=\delta_{ij}$, and let
$\Omega=(\Omega_{ij})_{i,j=1}^g$ be the matrix of $b$-periods:
\begin{equation*}
\Omega_{ij}:=\int_{b_i}\omega_j.
\end{equation*}
By Riemann bilinear relations, the matrix $\Omega$ is symmetric with positive definite
imaginary part
\begin{equation*}
\Imm\Omega:=\frac{\Omega-\bar{\Omega}}{2\sqrt{-1}}.
\end{equation*}
Denote by $\Lambda=\ZZ^g\oplus \Omega\ZZ^g$ the period lattice of a Torelli marked
Riemann surface $X$. It is a discrete subgroup of rank $2g$ in $\CC^g$. The Jacobian
variety of $X$ is defined as the following $g$-dimensional complex torus
\begin{equation*}
J(X):=\CC^g/\Lambda.
\end{equation*}
There is a complex-analytic isomorphism $J(X)\simeq\Alb(X)$ which is obtained by
choosing the basis $\{\int_{a_i}\}_{i=1}^g$ for $\check{H}^0(X,K)^\vee$ of the
$a$-periods. The complex coordinates $Z={}^t(z_1,\dotsc,z_g)$ on $\Alb(X)$ for this
basis are the standard complex coordinates on $\CC^g$ and $H_1(X,\ZZ)\simeq\Lambda$.
The invariant K\"{a}hler metric on $J(X)$ in these coordinates has the form:
\begin{equation*}
ds^2=\sum_{i,j=1}^g Y^{ij} dz_i\otimes d\bar{z}_j,
\end{equation*}
where $\{Y^{ij}\}_{i,j=1}^g=Y=(\Imm\Omega)^{-1}$, so that the corresponding Hermitian
inner product on $\CC^g=T_0 J(X)$ is given by
\begin{equation*}
<U,V>\,=\sum_{i,j=1}^g Y^{ij}u_i\bar{v}_j,
\end{equation*}
where $U={}^t(u_1,\dotsc,u_g), V={}^t(v_1,\dotsc,v_g)$.

The Abel-Jacobi map $\mu:\Div^0(X)\rightarrow J(X)$ for a Torelli marked Riemann
surface $X$ is defined by
\begin{equation*}
\mu(D)={}^t(\sum_{i=1}^n\int^{P_i}_{Q_i}\omega_1,\dotsc,
\sum_{i=1}^n\int^{P_i}_{Q_i}\omega_g)\in J(X),
\end{equation*}
where $D=\sum_{i=1}^n(P_i-Q_i)\in\Div^0(X)$. By Abel's theorem, the Abel-Jacobi map
establishes a complex-analytic isomorphism $\Div^0(X)/\PDiv(X)\simeq J(X)$.

Two other maps, the holonomy map
\begin{equation*}
hol: \Div^0(X)\rightarrow\Hom(\pi_1(X),U(1))=\Hom(H_1(X,\ZZ),U(1)),
\end{equation*}
and the inclusion map
\begin{equation*}
\imath: H^1(X,\RR)/H^1(X,\ZZ)\rightarrow J(X),
\end{equation*}
are defined as follows. For $D=\sum_{i=1}^l n_i P_i\in\Div^0(X)$ let $\omega_D$ be the
unique meromorphic differential on $X$ with only simple poles at $P_i$ with residues
$n_i$ and with pure imaginary periods:
\begin{equation*}
\int_c\omega_D\in \sqrt{-1}\,\RR~\text{for all}~[c]\in H_1(X,\ZZ).
\end{equation*}
The differential of the third kind $\omega_D$ can be also interpreted in terms of the
holomorphic line bundle $L=[D]$ associated with the divisor $D$ as a unique
connection-current $\omega_D$ with the property that its curvature-current
$\frac{\sqrt{-1}}{2\pi}\bar{\partial}\,\omega_D$ is Poincar\'{e} dual to the $0$-cycle
$D$. Set
\begin{equation*}
\Pi_a={}^t(\int_{a_1}\omega_D,\,\dotsc,\int_{a_g}\omega_D),~
\Pi_b={}^t(\int_{b_1}\omega_D,\,\dotsc,\int_{b_g}\omega_D)\in\RR^g,
\end{equation*}
and define the holonomy map $hol(D)=\rho\in\Hom(\pi_1(X), U(1))$ by
\begin{equation*}
\rho(a_i)=e^{2\pi\sqrt{-1}\,\Pi_{ai}},~\rho(b_i)=e^{2\pi\sqrt{-1}\,\Pi_{bi}},~i=1,
\dotsc,g.
\end{equation*}
To define the inclusion map $\imath$ consider the isomorphism
\begin{equation*}
H^1(X,\RR)/H^1(X,\ZZ)\simeq\Hom(\pi_1(X),U(1)),
\end{equation*}
given by the exponential of the period map
\begin{equation*}
H^1(X,\RR)/H^1(X,\ZZ)\ni[\omega]\mapsto e^{2\pi\sqrt{-1}\int_c\omega}:=
e^{2\pi\sqrt{-1}\Pi_c}\in\Hom(\pi_1(X),U(1)),
\end{equation*}
and set
\begin{equation*}
\imath([\omega])=-\Omega\Pi_a+\Pi_b\bmod\Lambda.
\end{equation*}
The reciprocity law between differentials of the first and third kinds~(see,
e.g.~\cite{Kra}) gives
\begin{equation*}
\mu(D)=\imath(hol(D)).
\end{equation*}

These relations between the Abel-Jacobi map, the holonomy map and the inclusion map
can be summarized in the following commutative diagram
\begin{equation*}
\begin{CD}
\Pic^0(X)@>{\simeq}>>\Alb(X)@>{\simeq}>> \Div^0(X)/\PDiv(X) \\ @VVholV @VV{\nu}V
@VV{\mu}V \\ \Hom(\pi_1(X),U(1))@>{\simeq}>>H^1(X,\RR)/H^1(X,\ZZ)@>{\imath}>> J(X)
\end{CD}
\end{equation*}
where the isomorphism $\nu$ is defined by the commutativity of the diagram.

Let $\theta(Z\,|\,\Omega)$ be the Riemann theta-function of a Torelli marked Riemann
surface $X$, defined by the following absolutely convergent series,
\begin{equation*}
\theta(Z\,|\,\Omega):=\sum_{n\in\ZZ^g}e^{\pi\sqrt{-1}((\Omega
n,n)+2(n,Z))},~Z\in\CC^g,
\end{equation*}
where $(A,B)={}^t A B$ for $A,B\in\CC^g$, and denote by $\theta[\xi](Z\,|\,\Omega)$
the theta-function with characteristics $\xi={}^t(\xi_a,\xi_b)\in\RR^{2g}$
\begin{equation*}
\theta[\xi](Z\,|\,\Omega):=e^{\pi\sqrt{-1}((\Omega\xi_a,\xi_a) +
2(\xi_a,Z+\xi_b))}\theta(Z+\Omega\xi_a+\xi_b\,|\,\Omega).
\end{equation*}
The Riemann theta-function satisfies the modular transformation formula
\begin{equation*}
\theta(-\Omega^{-1} Z\,|\,-\Omega^{-1})=\left(\det\left(\frac{\Omega}
{\sqrt{-1}}\right)\right)^{1/2} e^{\pi\sqrt{-1}(\Omega^{-1}Z,Z)} \theta(Z\,|\,\Omega).
\end{equation*}

Let $\Theta\subset J(X)$ be theta-divisor --- the zero locus of Riemann theta-function
on $J(X)$. It can be defined geometrically as a divisor of the unique (up to a
translation) holomorphic line bundle over $J(X)$ whose first Chern class is given by
the intersection form in $H_1(J(X),\ZZ)=\Lambda\simeq H_1(X,\ZZ)$. The theta-divisor
is even, $\Theta=-\Theta$, and depends on the marking. Denote by $\mathcal{E}=\{L\in
\Pic^{g-1}(X)\,|\,h^0(L)>0\}$ the so-called canonical theta-divisor. It follows from
the Riemann-Roch formula that $L\in\mathcal{E}$ if and only if $K\otimes L^{-1}\in
\mathcal{E}$. Next, choose a base point $P_0\in X$ and denote by
$\mu_{g-1}:\Div^{g-1}(X)\rightarrow J(X)$ the corresponding Abel map:
$\mu_{g-1}(D)=\mu(D-(g-1)P_0)$ for $D\in\Div^{g-1}(X)$. Let
$W_{g-1}=\mu_{g-1}(\mathcal{E})=\mu_{g-1}(X^{(g-1)})\subset\Pic^{g-1}(X)$, where
$X^{(g-1)}$ is $g-1$-fold symmetric product of $X$ and we have used the isomorphism
$\Pic^{g-1}(X)=\Div^{g-1}(X)/\PDiv(X)$. According to the Riemann theorem, there exists
$\kappa\in J(X)$ such that $\Theta=W_{g-1}+\kappa$. A vector of Riemann constants
$\kappa$ depends on a Torelli marking of the Riemann surface $X$ and on the choice of
a base-point, and has the property that there exists $\mathcal{K}\in\Pic^{g-1}(X)$
satisfying $\mathcal{K}\otimes\mathcal{K}=K$ such that
$\kappa=-\mu_{g-1}(\mathcal{K})$. Using the spin structure $\mathcal{K}$ we identify
$\Pic^{g-1}(X)\simeq \Pic^0(X)$ by $L\mapsto L\otimes\mathcal{K}^{-1}$.

\subsection{Green's functions}
For any Hermitian (i.e.~conformal) metric $ds^2$ on a Riemann surface $X$ let
\begin{equation*}
2\Delta_0=-\ast\bar{\partial}\ast \bar{\partial}
\end{equation*}
be the $\bar{\partial}$-Laplacian acting on functions on $X$. Here $\ast$ is the Hodge
star-operator for the metric $ds^2$, and $\bar{\partial}$ is the $(0,1)$-component of
the de Rham differential $d$ on $X$, $d=\partial + \bar{\partial}$. In local
coordinates $ds^2=\rho(z)|dz|^2$ and $\Delta_0=-\rho(z)^{-1}\partial^2/\partial
z\partial\bar{z}$. The Green's function $G$ of the Laplacian $\Delta_0$ is canonically
defined by the following properties~(see, e.g.~\cite{Lang}).
\begin{itemize}
\item[\textbf{G1.}] $G\in C^{\infty}(X\times X\setminus\Delta, \RR)$, where $\Delta$ is
the diagonal in $X\times X$.
\item[\textbf{G2.}] For every $P\in X$ there exists a neighborhood $U\ni P$ such
that the function
\begin{equation*}
\tilde{G}(P,Q):=G(P,Q)+\frac{1}{\pi}\log|z(P)-z(Q)|^2
\end{equation*}
is smooth in $U\times U$, where $z$ is a local coordinate at $U$.
\item[\textbf{G3.}] For every $Q\in X$ the function $G_Q(P):=G(P,Q)$ on $X\setminus\{Q\}$
satisfies
\begin{equation*}
\Delta_0 G_Q=-\frac{1}{\Area(X)},
\end{equation*}
where
\begin{equation*}
\Area(X)=\int_X\ast 1
\end{equation*}
is the area of $X$.
\item[\textbf{G4.}]
\begin{equation*}
\int_X \ast\,G_Q=0
\end{equation*}
for every $Q\in X$.
\end{itemize}
In distributional sense, conditions~\textbf{G1-G3} can be summarized in the following
single equation
\begin{equation*}
\Delta_0 G(P,Q)=\delta(P-Q)-\frac{1}{\Area(X)},
\end{equation*}
where Laplacian $\Delta_0$ acts on the first argument of $G$.

It follows from~\textbf{G1-G4} that the Green's function is symmetric:
$G(P,Q)=G(Q,P)$.

Next, consider the following tensor on $X\times X$
\begin{equation*}
S:=-\pi\partial\partial^\prime G=-\pi\frac{\partial^2 G}{\partial z\partial
w}dz\otimes dw,
\end{equation*}
where $\partial$ and $\partial^\prime$ act on the first and second arguments of $G$
correspondingly. It follows from \textbf{G1-G3} that $S$ is a symmetric bidifferential
of the second kind on $X\times X$ with biresidue 1. Specifically, $S$ is a symmetric
section of the line bundle $K_X \boxtimes K_X$ over $X\times X$, holomorphic on
$X\times X\setminus\Delta$ and having a pole of order $2$ at $\Delta$ with residue
$1$:
\begin{equation*}
S(P,Q)=\frac{dz(P)\otimes dz(Q)}{(z(P)-z(Q))^2}+O(1)~\text{as}~P\rightarrow Q.
\end{equation*}
Succinctly, $S\in\check{H}^0 (X\times X, K_X\boxtimes K_X(2\Delta))^{S_2}$ and
$\Bires|_\Delta S=1$. As it follows from~\textbf{G1-G3} and the Stokes theorem,
bidifferential $S$ has the property
\begin{equation*}
\text{v.p.}\int_X S_Q\,\bar\omega=0~\text{for all}~Q\in X~\text{and}~\omega
\in\check{H}^0(X,K).
\end{equation*}
Here $S_Q=\imath_Q^*S\in\check{H}^0(X,K(2Q))$ is a pull-back of $S$ by the map
$\imath_Q: X\hookrightarrow X\times X$ defined by $\imath_Q(P)=(P,Q)$, and the
integral is understood in the principal value sense. This property uniquely
characterizes $S$ as a classical Schiffer kernel. The Schiffer kernel is defined
canonically and does not depend on the choice of a conformal metric $ds^2$ on $X$. As
it follows from the definition of $S$ and~\textbf{G1-G3}, it can be also uniquely
characterized by the property that for all $Q\in X$ the differential of the second
kind $S_Q$ on $X$ has pure imaginary periods.

Another symmetric bidifferential $B$ of the second kind with biresidue $1$ is
canonically associated with a Torelli marked Riemann surface $X$ and is uniquely
characterized by the property that for all $Q\in X$ the differential of the second
kind $B_Q=\imath_Q^*B$ on $X$ has zero $a$-periods. Remaining $b$-periods of $B$
satisfy the relations
\begin{equation*}
\int_{b_i}B_Q=2\pi\sqrt{-1}\,\omega_i(Q),~i=1,\dotsc,g,
\end{equation*}
which follow from the reciprocity law between differentials of the first and second
kinds (see, e.g.~\cite{G-H,Kra}).

The Schiffer kernel and the $B$-kernel are related as follows~\cite{Fay2}.
\begin{lemma} (The Fay's formula)
\begin{equation*}
B(P,Q)=S(P,Q)+\pi\sum_{i,j=1}^g Y^{ij}\,\omega_i(P)\,\omega_j(Q).
\end{equation*}
\end{lemma}
\begin{proof} It is sufficent to verify that the right hand side of the Fay's formula
has zero $a$-periods. Writing $S=d^\prime (-\pi \partial G)+\tilde{S}$, where
$\tilde{S}:=\pi\partial\bar{\partial}^\prime G$ is a regular bidifferential on
$X\times X$ holomorphic with respect to $P$ and anti-holomorphic with respect to $Q$,
we get by Stokes theorem that the periods of $S$ with respect to the variable $Q$ are
the same as periods of $\tilde{S}$. Next, it follows from properties \textbf{G1-G3}
and the Stokes theorem that
\begin{align*}
\int_X \tilde{S}_P\omega & =-\pi\omega(P)~\text{for all}~P\in X~\text{and}~
\omega\in\check{H}^0(X,K), \\ \int_X \tilde{S}_P\partial f &=0~\text{for all}~P\in
X~\text{and}~f\in C^\infty(X).
\end{align*}
Here $\tilde{S}_P=\jmath_P^*\tilde{S}\in\check{H}^0(X,\overline{K})$ is a pull-back of
the bidifferential $\tilde{S}$ by the map $\jmath_P: X\hookrightarrow X\times X$,
defined by $\jmath_P(Q)=(P,Q)$. These properties uniquely characterize
$K:=-\tilde{S}/\pi$ as a kernel of the Hodge projection operator
$P:L_2^{1,0}(X)\rightarrow H^{1,0}_{\bar\partial}(X)$ onto the subspace of harmonic
$(1,0)$-forms on $X$, so that
\begin{equation*}
K(P,Q)=\sum_{i,j=1}^g Y^{ij}\omega_i(P)\bar{\omega}_j(Q).
\end{equation*}
Therefore
\begin{equation*}
\int_{a_j} \tilde{S}_P=-\pi\sum_{i=1}^g Y^{ij}\omega_i(P),~j=1,\dotsc,g,
\end{equation*}
and the $a$-periods of the right hand side of the Fay's formula are indeed $0$.
\end{proof}

\begin{remark}
Using classical terminology (see, e.g.~\cite{Kra}), the kernel $K$ is the Bergman
reproducing kernel for the space of $(1,0)$-forms on the Riemann surface $X$.
\end{remark}

\section{Bosonic action functional and partition function}

For a \emph{classical $U(1)$-valued bosonic field} $g\in C^\infty(X,U(1))$ define
corresponding \emph{field current} as the following $1$-form on $X$
\begin{equation*}
\J:=\frac{1}{2\pi\sqrt{-1}}g^{-1}d g.
\end{equation*}
The 1-form $\J$ is closed, real-valued, and has integral periods: $[\J]\in
H^1(X,\ZZ)$. Denoting by $\mathcal{J}(X)$ the set of all 1-forms on $X$ with these
properties, we have canonical isomorphism
\begin{equation*}
C^\infty(X,U(1))/U(1)\simeq\mathcal{J}(X),
\end{equation*}
where the inverse map is given by the integration:
\begin{equation*}
g(P)=e^{2\pi\sqrt{-1}\int_{P_0}^P\J}.
\end{equation*}
Consider the following functional on the space of field currents $\mathcal{J}(X)$
\begin{equation*}
S(\J):=\frac{1}{4}\int_X \J\wedge\ast \J,
\end{equation*}
where $\ast$ is the Hodge star-operator (which for 1-forms on $X$ does not depend on
the choice of a conformal metric on $X$). According to the Hodge decomposition, every
$\J\in\mathcal{J}(X)$ can be uniquely written as
\begin{equation*}
\J=d\varphi_0+h,
\end{equation*}
where $\varphi_0\in C^\infty(X,\RR)$ and $h\in \mathcal{H}^1(X,\ZZ)$ is harmonic
$1$-form, $dh=d\ast h=0$, with integral periods. Therefore
\begin{equation*}
\mathcal{J}(X)=C^\infty(X,\RR)/\RR \times \mathcal{H}^1(X,\ZZ),
\end{equation*}
so that
\begin{equation*}
S(\J)=\frac{1}{4}\int_X d \varphi_0\wedge\ast d\varphi_0 +\frac{1}{4}\int_X
h\wedge\ast h:=S_0(\varphi_0)+S_{inst}(h).
\end{equation*}
Here the term $S_0$ is interpreted as an \emph{action functional} of the standard
theory of \emph{free bosons} on $X$, and the term $S_{inst}$ --- as a contribution
from \emph{instantons}.
\begin{remark}
The field current $\J$ can be also written as $\J=d\varphi$, where $\varphi$ is an
additive multi-valued real function on $X$ with integral periods, i.e. $\varphi$ is a
single-valued function on the universal cover $\tilde{X}$ of $X$ satisfying
$\varphi\circ\gamma-\varphi\in\ZZ$ for all $\gamma\in\pi_1(X)$, where $\pi_1(X)$ acts
on $\tilde{X}$ by deck transformations. The classical field $\varphi$ with these
properties is the simplest example of an \emph{instanton configuration} with the
\emph{instanton numbers} given by the periods of the 1-form $\J=d\varphi$.
\end{remark}

According to the \emph{bosonization procedure}~\cite{AG-B-M-N-V,K-N-T-Y}, quantum
theories of free bosons on $X$ are parameterized by the set
$\Pic^{g-1}(X)\setminus\mathcal{E}$ of generalized spin structures without \emph{zero
modes}. For every $L\in\Pic^{g-1}(X)\setminus\mathcal{E}$,
following~\cite{AG-B-M-N-V}, define the \emph{topological term} of the bosonic action
functional by
\begin{equation*}
S_{top}(\J):=\sqrt{-1}\int_X \J\wedge\theta_L+\frac{\sqrt{-1}}{2}\varepsilon(h)
=\sqrt{-1}\int_X h\wedge\theta_L+\frac{\sqrt{-1}}{2}\varepsilon(h).
\end{equation*}
Here (see Section 2.1)
\begin{equation*}
[\theta_L]:=hol(L_0)\in H^1(X,\RR)/H^1(X,\ZZ)\simeq \Hom(\pi_1(X),U(1)),
\end{equation*}
where $L_0=L\otimes\mathcal{K}^{-1}$, and $\varepsilon(h)$ is a parity of the spin
structure associated with $\imath(h)/2$:
\begin{equation*}
\varepsilon(h):=(l,m)\bmod 2~\text{where}~\imath(h)=-\Omega
l+m\in\Lambda,~l,m\in\ZZ^g.
\end{equation*}
Since $h\in H^1(X,\ZZ)$, the exponential $\exp\{2\pi S_{top}(\J)\}$ does not depend on
the choice of representative $\theta_L$.

The total action functional of the theory of free bosons on $X$ with a topological
term parameterized by $L\in\Pic^{g-1}(X)\setminus\mathcal{E}$ is given by
\begin{equation*}
S_L(\J):=S(\J)+S_{top}(\J)=S_0(\varphi_0)+S_{inst}(h)+S_{top}(h),
\end{equation*}
and is, in general, complex-valued. The \emph{partition function} of the theory is
defined by the following functional integral
\begin{equation*}
\langle\mathbf{1}\rangle_L:=\int_{\mathcal{J}(X)}[\mathcal{D}\J]\,e^{-2\pi S_L(\J)}.
\end{equation*}

Mathematically rigorous definition is the following. Using the Hodge decomposition,
set
\begin{equation*}
\langle\mathbf{1}\rangle_L:=\mathbf{Z}_0 \mathbf{Z}_{inst}.
\end{equation*}
Here
\begin{equation*}
\mathbf{Z}_0:=\int_{C^\infty(X,\RR)/\RR}[\mathcal{D}\varphi_0] e^{-2\pi
S_0(\varphi_0)}
\end{equation*}
is a \emph{fluctuation part} of $\langle\mathbf{1}\rangle_L$
--- the partition function of the quantum field theory of free bosons on
$X$ given by the standard Gaussian integral, and
\begin{equation*}
\mathbf{Z}_{inst} :=\sum_{h\in\mathcal{H}^1(X,\ZZ)} e^{-\frac{\pi}{2}\int_X
h\wedge\ast h -2\pi \sqrt{-1}\int_X h\wedge\theta_L -\pi\sqrt{-1}\varepsilon(h)}
\end{equation*}
is an \emph{instanton part} of $\langle\mathbf{1}\rangle_L$ with topological term. The
instanton part $\mathbf{Z}_{inst}$ is well-defined: it is given by the absolutely
convergent series over the lattice $\mathcal{H}^1(X,\ZZ)$ of rank $2g$ in $\RR^{2g}$.
In order to define the Gaussian integral for $\mathbf{Z}_0$ choose a conformal metric
$ds^2$ on $X$ so that the functional $S_0$ becomes a quadratic form of the Laplacian
$\Delta_0$. According to~\cite{Polyakov}, metric $ds^2$ defines a ``Riemannian
metric'' on the infinite-dimensional Frech\'{e}t manifold $C^\infty(X,\RR)$
\begin{equation*}
||v||^2:=\int_X \ast\,|v|^2,~v\in T_{\varphi_0}C^\infty(X,\RR),
\end{equation*}
and the integration measure $[\mathcal{D}\varphi_0]$ on $C^\infty(X,\RR)$ is defined
as the volume form of this metric. Mathematically this is equivalent (see,
e.g.~\cite{Schvartz}) to a definition
\begin{equation*}
\mathbf{Z}_0=\int_{C^\infty(X,\RR)/\RR}[\mathcal{D}\varphi_0] e^{-2\pi\int_X
\Delta_0\varphi_0\ast\varphi_0}:=
\left(\frac{\Area(X)}{\det_{\zeta}\Delta_0}\right)^{1/2}.
\end{equation*}
Here $\det_{\zeta}\Delta_0$ is a functional determinant of the Laplacian $\Delta_0$ of
the metric $ds^2$ defined through the zeta-function of $\Delta_0$. The area term ---
contribution from \emph{zero modes} --- reflects the integration over the coset
$C^\infty(X,\RR)/\RR$. This completes the rigorous definition of the partition
function $\langle\mathbf{1}\rangle_L$.
\begin{remark} Though the functional $S_0$ does not depend on the choice of the
conformal metric $ds^2$, this metric is required for defining the integration measure
$[\mathcal{D}\varphi_0]$, and the determinant $\det_\zeta\Delta_0$ depends on $ds^2$.
This is Polyakov's \emph{conformal anomaly}~\cite{Pol2}, which is due to the fact that
any ``regularization procedure'' defining $\det\Delta_0$ breaks conformal invariance.
\end{remark}

It is remarkable that the instanton part $\mathbf{Z}_{inst}$ can be computed
explicitly in terms of the Riemann theta-function with characteristics
$\xi={}^t(-\Pi_a,\Pi_b)\in\RR^{2g}$, where
\begin{equation*}
e^{2\pi\sqrt{-1}\Pi_{ai}}:=hol(L_0)(a_i),~e^{2\pi\sqrt{-1}\Pi_{bi}}:=
hol(L_0)(b_i),~i=1,\dotsc,g.
\end{equation*}
Namely, the following statement holds~\cite{AG-B-M-N-V}.
\begin{proposition} \label{theta}
\begin{equation*}
\mathbf{Z}_{inst}= \frac{2^{g/2}}{\sqrt{\det
Y}}e^{\frac{\pi}{2}<Z-\bar{Z},Z-\bar{Z}>}|\theta(Z\,|\,\Omega)|^2
=\frac{2^{g/2}}{\sqrt{\det Y}}|\theta[\xi](0\,|\,\Omega)|^2,
\end{equation*}
where $Z:=\imath(\theta_L)=\Omega\xi_a+\xi_b=-\Omega\Pi_a+\Pi_b\in\CC^g$.
\end{proposition}
\begin{proof} It is a straightforward application~(cf.~\cite{AG-B-M-N-V}) of the Poisson
summation formula
\begin{equation*}
\sum_{n\in\ZZ^{2g}}f(n)=\sum_{n\in\ZZ^{2g}}\hat{f}(n),
\end{equation*}
where $f\in\mathcal{S}(\RR^{2g})$ is a function of the Schwartz class and
\begin{equation*}
\hat{f}(p):=\int_{\RR^{2g}}e^{-2\pi\sqrt{-1}(x,p)}f(x)d^{2g}x
\end{equation*}
is its Fourier transform, where $(~,~)$ is the standard Euclidean inner product in
$\RR^{2g}$. Let $\{\alpha_i,\beta_i\}_{i=1}^{g}$ be the basis for the lattice
$\mathcal{H}^1(X,\ZZ)$ dual to a symplectic basis $\{a_i,b_i\}_{i=1}^{g}$ in
$H_1(X,\ZZ)$, so that for $h\in\mathcal{H}^1(X,\ZZ)$
\begin{equation*}
h=\sum_{i=1}^g(l_i\alpha_i + m_i\beta_i)~\text{and}~
\theta_L=\sum_{i=1}^g(\xi_{ai}\alpha_i + \xi_{bi}\beta_i).
\end{equation*}
Using $\varepsilon(h)=(l,m)\bmod 2$ and the formulas
\begin{equation*}
\int_X h\wedge\ast h=<\lambda,\lambda>,~\,2\sqrt{-1}\int_X h\wedge\theta_L=
<Z,\lambda>-<\lambda,Z>,
\end{equation*}
where $\lambda:=-\Omega l+m\in\Lambda$, the instanton part $\mathbf{Z}_{inst}$ can be
represented as the following theta-series
\begin{equation} \label{theta-series}
\mathbf{Z}_{inst}=\sum_{n\in\ZZ^{2g}}e^{-\frac{\pi}{2}(Q n,n) - \pi(A,n)}.
\end{equation}
Here
\begin{equation*}
n=
\begin{pmatrix}
  l \\
  m
\end{pmatrix}
 \in\ZZ^{2g},
\end{equation*}
$Q$ is the following $2g\times 2g$ matrix
\begin{equation*}
Q=
\begin{pmatrix}
  \bar{\Omega}Y\Omega & -\bar{\Omega}Y \\
  -Y\Omega & Y
\end{pmatrix},
\end{equation*}
and
\begin{equation*}
A=
\begin{pmatrix}
  \Omega Y\bar{Z}-\bar\Omega Y Z \\
  Y(Z-\bar{Z})
\end{pmatrix}
\in\CC^{2g}.
\end{equation*}
We apply the Poisson summation formula to the function
\begin{equation*}
f(x)=e^{-\frac{\pi}{2}(Q x,x) -\pi(A,x)}.
\end{equation*}
The inverse matrix $Q^{-1}$ is readily computed
\begin{equation*}
Q^{-1}=\frac{\sqrt{-1}}{2}
\begin{pmatrix}
  \Omega^{-1} & I \\
  I & \bar{\Omega}
\end{pmatrix}
\end{equation*}
(this is the place where the extra term $\varepsilon(h)$ in the topological action is
crucial, cf.~\cite{AG-B-M-N-V}), so that
\begin{equation*}
\det Q=2^g \det Y\det(-\sqrt{-1}\Omega).
\end{equation*}
The Gaussian integral for $\hat{f}$ is computed explicitly
\begin{equation*}
\hat{f}(x)=\biggl(\frac{2^{g}}{\det(-\sqrt{-1}Y\Omega)}\biggr)^{1/2}e^{
\frac{\pi}{2}(Q^{-1} A,A) -2\pi(Q^{-1}(x+\sqrt{-1}A),x)},
\end{equation*}
so that
\begin{align*}
\sum_{n\in\ZZ^{2g}}\hat{f}(n)&=\biggl(\frac{2^{g}}{\det
(-\sqrt{-1}Y\Omega)}\biggr)^{1/2}
e^{\frac{\pi}{2}(<Z-\bar{Z},Z-\bar{Z}>-2\sqrt{-1}(\Omega^{-1}Z,Z))}\\
&\bar{\theta}(Z|\Omega)\theta(-\Omega^{-1}Z\,|\,-\Omega^{-1}).
\end{align*}
Using the modular transformation formula for the Riemann theta-function completes the
proof.
\end{proof}
\begin{remark} Note that though the action functional $S_L(\J)$ is not real-valued,
Proposition 2.3 shows that the partition function $\langle\mathbf{1}\rangle_L$ is real
and positive. Since under the involution $L\mapsto K\otimes L^{-1}$, or equivalently,
$Z\mapsto -Z$, $S_L(\J)\mapsto\overline{S_L(\J)}$ --- the complex conjugate, the
partition function $\langle\mathbf{1}\rangle_L$ has a $\ZZ_2$-symmetry $Z\mapsto -Z$,
in accordance with $\theta(Z\,|\,\Omega)=\theta(-Z\,|\,\Omega)$.
\end{remark}
\begin{remark}
As it follows from the Proposition 2.3 and~~\cite{AG-B-M-N-V,Fay1}, the partition
function $\langle\mathbf{1}\rangle_L$ coincides with the analytic torsion $T(L)$ for
the unitary line bundle $L_0=L\otimes\mathcal{K}^{-1}$ --- the zeta-function
determinant of the $\bar\partial$-Laplacian acting on sections of the line bundle $L$.
This is the essence of the spin-$1/2$ bosonization formula~\cite{AG-M-V,B-N}.
\end{remark}

Let $\tau(\mathbf{t},\mathbf{U}),\,\mathbf{t}=(t_1,t_2,\dotsc),$ be Sato's
$\tau$-function on the universal Grassmannian manifold, defined as a solution of
Hirota's bilinear equations~(see, e.g.~\cite{K-N-T-Y}). It is
well-known~\cite{AG-G-R,V,K-N-T-Y} that if a point $\mathbf{U}$ in the universal
Grassmannian manifold is the image of a Torelli marked Riemann surface $X$ with the
line bundle $L\in\Pic^{g-1}(X)\setminus \mathcal{E}$ under the Krichever map, then
\begin{equation*}
\tau(\mathbf{t},\mathbf{U})|_{\mathbf{t}=0}:=\tau(Z)=\theta[\xi](0\,|\,\Omega),
\end{equation*}
$Z=\Omega\xi_a + \xi_b$.

We define the $\tau$-function $\btau(Z)$ of a Torelli marked Riemann surface $X$ with
the line bundle $L\in\Pic^{g-1}(X)\setminus \mathcal{E}$ as the partition function
$\langle\mathbf{1}\rangle_L$. According to the Proposition 2.3,
\begin{equation*}
\btau(Z):=\langle\mathbf{1}\rangle_L=\frac{2^{g/2}} {\sqrt{\det Y}}|\tau(Z)|^2.
\end{equation*}

\section{$U(1)$-gauge symmetry Ward identities}

Introduce holomorphic and anti-holomorphic components of the bosonic field current
$\J$ through the decomposition
\begin{equation*}
\J=\frac{\sqrt{-1}}{2\pi}(J + \bar{J}),
\end{equation*}
so that $J=-g^{-1}\partial g~\text{and}~\bar{J}=-g^{-1}\bar{\partial} g$. By
definition, \emph{multi-point correlation functions} of current components
$J~\text{and}~\bar{J}$ are given by the following functional integral
\begin{multline*}
\langle J(P_1)\dotsb J(P_m)\bar{J}(Q_1)\dotsb\bar{J}(Q_n)\rangle \\
:=\int_{\mathcal{J}(X)}[\mathcal{D}\J] \,J(P_1)\dotsb
J(P_m)\bar{J}(Q_1)\dotsb\bar{J}(Q_n)\,e^{-2\pi S(\J)}.
\end{multline*}

We start with the computation of the \emph{normalized} $1$-point correlation function
of the holomorphic component $J$, defined by
\begin{equation*}
\langle\langle J(P)\rangle\rangle:=\frac{\langle J(P)\rangle}{\langle\mathbf{1}
\rangle_L}.
\end{equation*}
First observe that for the standard quantum field theory of free bosons on $X$ with
the action functional $S_0(\varphi_0)$ all multi-point current correlation functions
with odd number of components $\partial\varphi_0~\text{and}~\bar\partial\varphi_0$
vanish. This is a standard fact from the theory of Gaussian integration (for
mathematical treatment see, e.g.~\cite{G-J}). In particular,
$<\partial\varphi_0(P)>=0$ and from the Hodge decomposition
$J=-2\pi\sqrt{-1}(\partial\varphi_0+h^{1,0})$ it follows that
\begin{equation*}
\langle\langle J(P)\rangle\rangle=-2\pi\sqrt{-1}\langle\langle
h^{1,0}(P)\rangle\rangle.
\end{equation*}
Second, consider the partial derivative $\partial\btau/\partial z_i$ of the
$\tau$-function. In this case differentiation under the functional integral sign is
legitimate since the instanton part $\mathbf{Z}_{inst}$ --- the only part that
actually depends on $Z$ --- is given by the absolutely convergent series. Using the
explicit form of the topological term $S_{top}(\J)$ and the equation
\begin{equation*}
-\frac{\partial<Z,\lambda>}{\partial z_i}=-\sum_{j=1}^gY^{ij}\bar{\lambda}_j
=-2\sqrt{-1}\int_{a_i}h^{1,0},
\end{equation*}
we obtain
\begin{equation*}
\frac{\partial\log\btau}{\partial z_i}= \int_{a_i}\langle\langle J(P)\rangle\rangle,
~i=1,\dotsc,g.
\end{equation*}
This is a \emph{$U(1)$-gauge symmetry Ward identity} for the 1-point correlation
function $\langle\langle J(P)\rangle\rangle$. Using the formula
\begin{equation*}
\frac{\partial\log\btau}{\partial z_i}=\sum_{j=1}^g\pi Y^{ij}(z_j-\bar{z}_j)+
\frac{\partial\log\theta(Z\,|\,\Omega)}{\partial
z_i}=\frac{\partial\log\theta[\xi](U\,|\,\Omega)}{\partial u_i}\bigg|_{U=0},
\end{equation*}
where $U={}^t(u_1,\dotsc, u_g)$, and taking into account that $\langle\langle
J(P)\rangle\rangle$ is an abelian differential on $X$, we finally obtain
\begin{equation*}
\langle\langle J(P)\rangle\rangle=\sum_{i=1}^g
\frac{\partial\log\theta[\xi](U\,|\,\Omega)} {\partial u_i}\bigg|_{U=0}\omega_i(P).
\end{equation*}
This expression is in perfect agreement with spin-$1/2$ bosonization formulas
in~\cite{R2,K-N-T-Y}.

\begin{remark} In this derivation it was assumed that $1$-point correlation function
$\langle\langle J(P)\rangle\rangle$ is holomorphic on $X$. It is one of the basic
principles of conformal field theory that multi-point correlation functions of
holomorphic and anti-holomorphic current components are multi-differentials on $X$,
holomorphic and anti-holomorphic with respect to the corresponding variables. In our
case the action functional $S(\J)$ is quadratic and the statement immediately follows
from the representation for $\langle h^{1,0}(P)\rangle$, since $h$ is harmonic.
\end{remark}

Similarly,
\begin{equation*}
\frac{\partial<\lambda,Z>}{\partial \bar{z}_i}=\sum_{j=1}^gY^{ij}\lambda_j
=-2\sqrt{-1}\int_{a_i}h^{0,1},
\end{equation*}
and we obtain
\begin{equation*}
\frac{\partial\log\btau}{\partial \bar{z}_i}= \int_{a_i}\langle\langle
\bar{J}(P)\rangle\rangle,~i=1,\dotsc,g
\end{equation*}
--- a \emph{$U(1)$-gauge symmetry Ward identity} for the 1-point correlation
function $\langle\langle \bar{J}(P)\rangle\rangle$. As before,
\begin{equation*}
\langle\langle \bar{J}(P)\rangle\rangle=\sum_{i=1}^g
\frac{\partial\log\overline{\theta[\xi](U\,|\,\Omega)}} {\partial \bar{u}_i}
\bigg|_{U=0}\bar{\omega}_i(P).
\end{equation*}
\begin{remark}
There is no contradiction between the equations $\bar{J}(P)=-\overline{J(P)}$ and
$\langle\langle \bar{J}(P)\rangle\rangle = \overline{\langle\langle
J(P)\rangle\rangle}$ because the action functional $S_L(\J)$ is not real-valued.
Indeed, introducing explicit $Z$-dependence into the action functional and correlation
functions, we get
\begin{align*}
\langle\langle \bar{J}(P)\rangle\rangle_Z & =\int_{\mathcal{J}(X)}[\mathcal{D}\J]
\,\bar{J}(P)\,e^{-2\pi S(\J;Z)} \\ & =-\int_{\mathcal{J}(X)}[\mathcal{D}\J]
\,\overline{J(P)}\,e^{-2\pi \overline{S(\J;-Z)}}\\& =-\overline{\langle\langle
J(P)\rangle\rangle_{-Z}}= \overline{\langle\langle J(P)\rangle\rangle_{Z}},
\end{align*}
since $\partial\log\theta(Z\,|\,\Omega)/\partial z_i$ is odd.
\end{remark}

Now consider the normalized \emph{reduced} 2-point correlation function
\begin{equation*}
\langle\langle J(P)J(Q)\rangle\rangle:= \frac{\langle
J(P)J(Q)\rangle}{\langle\mathbf{1}\rangle_L} -\langle\langle J(P)\rangle\rangle
\langle\langle J(Q)\rangle\rangle.
\end{equation*}
First observe that
\begin{equation*}
\langle\langle J(P)J(Q)\rangle\rangle =
-4\pi^2(\langle\langle\partial\varphi_0(P)\partial\varphi_0(Q)\rangle\rangle +
\langle\langle h^{1,0}(P)h^{1,0}(Q)\rangle\rangle).
\end{equation*}
Next, use another standard fact from the theory of Gaussian integration~\cite{G-J}:
\begin{equation*}
\langle\langle\partial\varphi_0(P)\partial\varphi_0(Q)\rangle\rangle=
\frac{1}{4\pi}\partial\partial^\prime G(P,Q)=-\frac{1}{4 \pi^2}S(P,Q),
\end{equation*}
and consider the second partial derivative $\partial^2\log\btau/\partial z_i\partial
z_j$. Arguing as before we obtain
\begin{equation*}
\frac{\partial^2\log\btau}{\partial z_i\partial z_j}=
\int_{a_i}\int_{a_j}(\langle\langle J(P)J(Q)\rangle\rangle-S(P,Q)),~i,j=1,\dotsc,g.
\end{equation*}
This is a $U(1)$-gauge symmetry Ward identity for the $2$-point correlation function
$\langle\langle J(P)J(Q)\rangle\rangle$. It follows from the derivation and discussion
in Section 2.2 that the integrand is a holomorphic symmetric bidifferential on
$X\times X$. Therefore $\langle\langle J(P)J(Q)\rangle\rangle$ is a symmetric
bidifferential of the second kind on $X\times X$ with biresidue 1 on the diagonal.
Using the equation
\begin{equation*}
\frac{\partial^2\log\btau}{\partial z_i\partial z_j}=\pi Y^{ij}+
\frac{\partial^2\log\theta(Z\,|\,\Omega)}{\partial z_i\partial z_j}= \pi Y^{ij}+
\frac{\partial^2\log\theta[\xi](U\,|\,\Omega)} {\partial u_i\partial u_j}\bigg|_{U=0}
\end{equation*}
and the Fay's formula, we finally obtain
\begin{equation*}
\langle\langle J(P)J(Q)\rangle\rangle
=B(P,Q)+\sum_{i,j=1}^g\frac{\partial^2\log\theta[\xi](U\,|\,\Omega)} {\partial
u_i\partial u_j}\bigg|_{U=0}\omega_i(P)\omega_j(Q).
\end{equation*}
Again, this expression is in perfect agreement with formulas in~\cite{R2,K-N-T-Y}.

Next, consider $2$-point correlation function $\langle\langle
J(P)\bar{J}(Q)\rangle\rangle$ with holomorphic and anti-holomorphic components.
Arguing as before we obtain
\begin{equation}
\frac{\partial^2\log\btau}{\partial z_i\partial \bar{z}_j}=
\int_{a_i}\int_{a_j}(\langle\langle J(P)\bar{J}(Q)\rangle\rangle-K(P,Q)),
\end{equation}
where
\begin{equation*}
\langle\langle\partial\varphi_0(P)\bar\partial\varphi_0(Q)\rangle\rangle=
\frac{1}{4\pi}\partial\bar{\partial}^\prime G(P,Q)=\frac{1}{4 \pi^2}K(P,Q),
\end{equation*}
and the Bergman reproducing kernel $K$ was introduced in Section 2.2. Using explicit
representation for $K$ and the formula
\begin{equation*}
\frac{\partial^2\log\btau}{\partial z_i\partial \bar{z}_j}=-\pi Y^{ij}
\end{equation*}
we get
\begin{equation*}
\langle\langle J(P)\bar{J}(Q)\rangle\rangle =0.
\end{equation*}

To summarize, introduce multi-point reduced normalized correlation functions, defined
inductively as follows~(see, e.g.~\cite{Tak})
\begin{multline*}
\langle\langle J(P_1)\dotsb J(P_m)\bar{J}(Q_1)\dotsb\bar{J}(Q_n)\rangle\rangle:= \\
\frac{\langle J(P_1)\dotsb J(P_m)\bar{J}(Q_1)\dotsb\bar{J}(Q_n)\rangle}{\langle
\mathbf{1}\rangle_L} \\ -\sum_{l=2}^{m+n}\sum_{I=I_1\bigcup \dotsb\bigcup I_l}
\langle\langle J(I_1)\rangle\rangle\dotsb\langle\langle J(I_l)\rangle\rangle.
\end{multline*}
Here $I=\{P_1,\dots,P_m;Q_1,\dotsc,Q_n\}$ and summation goes over all partitions of
the set $I$ into the non-empty subsets $I_1,\dotsc, I_l$; for every subset
$I_k=\{P_{i_1},\dots,P_{i_{m_k}}; Q_{j_1},\dotsc,Q_{j_{n_k}}\}$ of $I$
\begin{equation*}
\langle\langle J(I_k)\rangle\rangle:=\langle\langle J(P_{i_1})\dotsb J(P_{i_{m_k}})
\bar{J}(Q_{j_1})\dotsb\bar{J}(Q_{j_{n_k}}) \rangle\rangle.
\end{equation*}

We have the following result (cf.~\cite{K-N-T-Y}).
\begin{theorem}
\begin{multline*}
\langle\langle J(P_1)\dotsb J(P_m)\bar{J}(Q_1)\dotsb\bar{J}(Q_n)\rangle\rangle \\
=\sum_{i_1=1}^g\dotsb\sum_{i_m=1}^g\sum_{j_1=1}^g\dotsb\sum_{j_n=1}^g
\frac{\partial^{m+n}\log |\theta(U\,|\,\Omega)|^2} {\partial u_{i_1}\dotsb\partial
u_{i_m}\partial \bar{u}_{j_1}\dotsb\partial \bar{u}_{j_n}} \bigg|_{U=0}\\
\omega_{i_1}(P_1)\dotsb\omega_{i_m}(P_m)\bar{\omega}_{j_1}(Q_1)\dotsb
\bar{\omega}_{j_n}(Q_n) \\ +\delta_{m,2}\delta_{n,0}B(P_1,P_2) +
\delta_{m,0}\delta_{n,2}\overline{B(Q_1,Q_2)}.
\end{multline*}
\end{theorem}
\begin{proof} Directly follows from the definition of normalized reduced
correlations functions and the arguments used above. Note that when $m n>0$ all
reduced correlation functions vanish.
\end{proof}

\part{Free Bosons and Tau-Functions for Closed Smooth Jordan Curves}
\section{Mathematical set-up} Let $C$ be a smooth closed Jordan curve in the
complex plane $\CC$ encircling the origin $0$. Denote by $\Omega$ the interior of a
contour $C$ --- a simply connected domain in $\CC$ containing $0$ and bounded by $C$,
and denote by $\mathcal{C}$ the set of all such contours $C$.
\subsection{Riemann mapping theorem} Let $D=\{w\in\CC\,\,|\,|w|<1\}$ be the unit
disk in the complex plane $\CC$. By Riemann mapping theorem, there exists a conformal
isomorphism $F:\Omega\rightarrow D$, uniquely characterized by the conditions
$F(0)=0~\text{and}~F^\prime(0)>0$. By Carath\'{e}odory's correspondence of the
boundaries principle $F$ extends to the regular map between the closure $\bar{\Omega}$
of the domain $\Omega$ and the closed unit disk $\bar{D}$, and $F|_C$ is a $C^\infty$
- isomorphism between the contour $C=\partial\Omega$ and the unit circle $S^1=\partial
D$. The inverse map $f=F^{-1}:D\rightarrow\Omega$ is a univalent function on $D$,
smooth up to the boundary. The value $r=F^\prime(0)=f^\prime(0)$ is called the
conformal radius of $C$ with respect to $0$.

Denote by $\mathcal{C}_1$ the set of all contours of conformal radius $1$ and by
$\mathcal{F}_1$ --- the set of all univalent functions on $D$ which are smooth up to
the boundary and normalized by the conditions $f(0)=0~\text{and}~f^\prime(0)=1$. Let
$\mathcal{A}_1$ be the affine space of all holomorphic functions on $D$ which are
smooth up to the boundary and have the same normalization at $0$ as functions in
$\mathcal{F}_1$. The space $\mathcal{A}_1$ has a structure of an infinite-dimensional
complex manifold with Frech\'{e}t topology given by the uniform convergence with all
derivatives in the closed unit disk $\bar{D}$. Taylor coefficients of the power series
expansion at $z=0$ are natural global coordinates on $\mathcal{A}_1$. The set
$\mathcal{F}_1$ is open in $\mathcal{A}_1$ and, therefore, is an infinite-dimensional
complex manifold. By Riemann mapping theorem there is a bijection
\begin{equation*}
\mathcal{C}_1\simeq\mathcal{F}_1,
\end{equation*}
which turns $\mathcal{C}_1$ into an infinite-dimensional complex manifold as well.

Let $\Diff_+(S^1)$ be a group of all orientation preserving diffeomorphisms of the
unit circle $S^1$, with $S^1$ interpreted as a rotation subgroup of $\Diff_+(S^1)$. It
was observed by A.A.~Kirillov~\cite{Kirillov1} that there is a remarkable bijection
between $\mathcal{C}_1$ and the infinite-dimensional homogeneous manifold
$\Diff_+(S^1)/S^1$. It turns $\Diff_+(S^1)/S^1$ into an infinite-dimensional complex
manifold and gives rise to canonical isomorphisms
\begin{equation*}
\mathcal{C}_1\simeq\mathcal{F}_1\simeq \Diff_+(S^1)/S^1.
\end{equation*}

In addition to the conformal map $F$, Kirillov's construction of the correspondence
$C\mapsto\gamma$ uses a conformal map $G$ of the exterior domain $\CC\setminus\Omega$
of the contour $C$ onto the exterior domain $\CC\setminus D$ of the unit circle $S^1$,
normalized by the conditions
$G(\infty)=\infty~\text{and}~G^\prime(\infty)>0$\footnote{ The value
$d=1/G^\prime(\infty)$ is called the transfinite diameter of $C$.}. Namely, if
$g=G^{-1}$ is the inverse map to $G$, then $\gamma:=F\circ g\in\Diff_+(S^1)$. To prove
that $C\mapsto\gamma$ is onto and to construct the inverse map,
A.A.~Kirillov~\cite{Kirillov1} used the Riemann theorem that all complex structures on
the 2-dimensional sphere $S^2$ are isomorphic to the standard complex structure of the
Riemann sphere $\PP^1$.

The manifold $\Diff_+(S^1)/S^1$ admits an interpretation as a set of all complex
structures on the space $\Omega\,\RR^d$ of based loops in $\RR^d$ and as such plays a
fundamental role in the string theory~\cite{Bow1,Bow2}.

\begin{remark} Similarly~\cite{Kirillov-Yuriev,N-V} it can be shown that another
homogeneous space $\Diff_+(S^1)/\Mob(S^1)$ is also an infinite-dimensional complex
manifold and $\Diff_+(S^1)/S^1$ is a holomorphic disk fiber space over
$\Diff_+(S^1)/\Mob(S^1)$.

There is a natural inclusion
\begin{equation*}
\Diff_+(S^1)/\Mob(S^1)\hookrightarrow T(1),
\end{equation*}
where $T(1)$ is a classical Bers universal Teichm\"{u}ller space
\begin{equation*}
T(1):=\Homeo_{qs}(S^1)/\Mob(S^1)
\end{equation*}
(see~\cite{Bers} for definitions and details). The space $T(1)$ has a natural
structure of an infinite-dimensional complex manifold and it was shown by S.~Nag and
A.~Verjovsky~\cite{N-V} that the inclusion map is holomorphic (moreover,
$\Diff_+(S^1)/\Mob(S^1)$ is one leaf of a holomorphic foliation on $T(1)$). It was
also shown in~\cite{N-V} that the pull-back by the inclusion map of the canonical
Weil-Petersson metric on $T(1)$ to $\Diff_+(S^1)/\Mob(S^1)$ coincides with the
K\"{a}hler metric introduced by A.A.~Kirillov and D.V.
Yuriev~\cite{Kirillov1,Kirillov-Yuriev} via the orbit method.
\end{remark}

Let $\PP^1_C$ be the double of the open Riemann surface $\PP^1\setminus\Omega$,
defined by gluing together $\PP^1\setminus\Omega$ and its copy
$\overline{\PP^1\setminus\Omega}$ with the opposite complex structure along their
common boundary $C$. As a smooth manifold $\PP^1_C$ is diffeomorphic to the
two-dimensional sphere $S^2$. The complex structure on $\PP^1_C$ is defined by the
structure sheaf $\mathcal{O}=\{\mathcal{O}_U\}$ --- the sheaf of germs of holomorphic
functions on $\PP^1_C$. Namely, for any connected open set $U\subset\PP^1_C$ define
the complex vector space $\mathcal{O}_U$ as follows.
\begin{itemize}
\item[(a)] If
$U\subset\PP^1\setminus\Omega$ then $\mathcal{O}_U$ is a complex vector space of all
holomorphic functions on $U$ with respect to the standard complex structure on
$\PP^1$.
\item[(b)] If $U\subset\overline{\PP^1\setminus\Omega}$ then $\mathcal{O}_U$ is a
complex vector space of all anti-holomorphic functions on $U$ with respect to the
standard complex structure on $\PP^1$.
\item[(c)] If $U\cap C\neq\emptyset$ then, setting
$U_+=U\cap(\PP^1\setminus\Omega),\,U_-=U\cap(\overline{\PP^1\setminus\Omega})$, define
\begin{equation*}
\mathcal{O}_U=\{(f_+ ,f_{-})\,|\,f_{\pm}\in \mathcal{O}_{U_{\pm}}~\text{and}~
f_+(z)=f_-(\bar{z})~\text{for}~z\in U\cap C\}.
\end{equation*}
\end{itemize}
Finally, for any open set $U\subset\PP^1_C$ define $\mathcal{O}_U$ as a direct sum of
the corresponding complex vector spaces for the connected components of $U$. Clearly,
complex vector spaces $\{\mathcal{O}_U\}$ form a sheaf $\mathcal{O}$ which canonically
defines the complex structure on $\PP^1_C$. By the Riemann theorem $\PP^1_C$ is
complex-isomorphic to the Riemann sphere $\PP^1$, i.e.~there exists a global
meromorphic coordinate $\zeta$ on $\PP^1_C$ such that
$\zeta(\infty)=\infty,~\zeta(\overline{\infty})=0$. In terms of the conformal mapping
$G$ this coordinate is given explicitly by
\begin{equation*}
\zeta(z)=
\begin{cases} G(z) & \text{if}~
z\in\PP^1\setminus\Omega, \\ 1/\overline{G(z)} &
\text{if}~z\in\overline{\PP^1\setminus\Omega}.
\end{cases}
\end{equation*}

\subsection{Green's functions} Let $G$ be the Green's
function\footnote{It should be always clear from the content whether $G$ is a
conformal map --- a function of one variable, or the Green's function --- a function
of two variables.} of the $\bar\partial$-Laplacian $\Delta_0$ of the conformal metric
$ds^2$ on $\PP^1$ acting on functions on $\PP^1$, and let $S$ and $B$ be corresponding
bidifferentials (see Section 1.2 in Part 1). Since $\PP^1$ has genus 0, the Fay's
formula reads $S=B$ and we have explicitly
\begin{equation*}
S(z,w):=-\pi\frac{\partial^2 G(z,w)}{\partial z\partial w}dz\otimes dw=\frac{dz\otimes
dw}{(z-w)^2},~z,w\in\CC.
\end{equation*}
Similarly, since the subspace of harmonic $(1,0)$-forms on $\PP^1$ is $\{0\}$, the
Bergman reproduction kernel $K$ for $\PP^1$ vanishes
\begin{equation*}
K(z,\bar{w}):=-\frac{\partial^2 G(z,w)}{\partial z\partial\bar{w}}dz\otimes
d\bar{w}=0.
\end{equation*}

Let $G_{DBC}$ be the Green's function for the Laplacian $\Delta_0$ on the exterior
domain $\PP^1\setminus\Omega$ with the Dirichlet boundary condition $\varphi|_C=0$.
One has
\begin{equation*}
G_{DBC}(z,w)=\frac{2}{\pi}G_{cl}(z,w),
\end{equation*}
where $G_{cl}$ is the classical Green's function for the domain
$\PP^1\setminus\Omega$, uniquely characterized by the following properties (see
e.g.~\cite{Hille}).
\begin{itemize}
\item[\textbf{CG1.}] For every $w\in\CC\setminus\Omega$ function $G_{cl}(z,w) +
\log|z-w|$ is harmonic in $(\CC\setminus\Omega)\setminus\{w\}$.
\item[\textbf{CG2.}] $G_{cl}(z,w)=0$ for all $w\in\PP^1\setminus\Omega~\text{and}~z\in C$.
\item[\textbf{CG3.}] $G_{cl}(z,w)=\log|w|+V+O(1/z)$ as $z\rightarrow\infty$,
where $V$ is the so-called Robin's constant of the domain $\Omega$.
\end{itemize}
The classical Green's function for the domain $\PP^1\setminus\Omega$ can be written
explicitly as a pull-back by the conformal map $G$ of the classical Green's function
for the domain $\PP^1\setminus D$,
\begin{equation*}
G_{cl}(z,w)=\log\bigg|\frac{1-G(z)\overline{G(w)}}{G(z)-G(w)}\bigg|.
\end{equation*}

According to Section 1.2 in part 1, the Schiffer kernel for $\PP^1\setminus\Omega$ is
defined as a symmetric bidifferential $S$ of the second kind on $(\PP^1\setminus\Omega
)\times (\PP^1\setminus\Omega)$ with the property
\begin{equation*}
\text{v.p.}\int_{\PP^1\setminus\Omega}S(z,w)\wedge \overline{u(w)}=0~ \text{for
all}~z\in\PP^1\setminus\Omega~\text{and}~u\in\mathcal{H}^{1,0}(\PP^1\setminus\Omega),
\end{equation*}
where $\mathcal{H}^{1,0}(\PP^1\setminus\Omega)$ is the complex vector space of
holomorphic $(1,0)$-forms on the domain $\PP^1\setminus\Omega$. It terms of the
Green's function the Schiffer kernel is given by
\begin{equation*}
S(z,w)=-\pi \frac{\partial^2 G_{DBC}(z,w)}{\partial z\partial w}dz\otimes dw,
\end{equation*}
and in terms of the conformal map $G$ it has the form
\begin{equation*}
S(z,w)=\frac{G^\prime(z)G^\prime(w)}{(G(z)-G(w))^2}dz\otimes dw.
\end{equation*}

Correspondingly, the Bergman reproducing kernel $K$ for the domain
$\PP^1\setminus\Omega$ is defined as a regular bidifferential on
$(\PP^1\setminus\Omega)\times (\PP^1\setminus\Omega)$, holomorphic with respect to the
first variable and anti-holomorphic with respect to the second variable, satisfying
the property\footnote{In this part we revert to the standard complex analyst's
notation $i=\sqrt{-1}$.}
\begin{equation*}
\frac{1}{2i}\int_{\PP^1\setminus\Omega}K(z,\bar{w})\wedge u(w)=u(z)~\text{for
all}~z\in\PP^1\setminus\Omega~\text{and}~u\in\mathcal{H}^{1,0}(\PP^1\setminus\Omega).
\end{equation*}

The Bergman reproducing kernel $K$ can be also characterized as a kernel of the
projection operator from $L_2^{1,0}(\PP^1\setminus\Omega)$ onto the subspace
$\mathcal{H}^{1,0}(\PP^1\setminus\Omega)$, where the Hilbert space
$L_2^{1,0}(\PP^1\setminus\Omega)$ is the space of $(1,0)$-forms on
$\PP^1\setminus\Omega$ with the $L_2$-norm
\begin{equation*}
||u||^2:=\frac{1}{2}\int_{\PP^1\setminus\Omega}u\wedge \ast u=
\frac{i}{2}\int_{\PP^1\setminus\Omega}u\wedge \bar{u}.
\end{equation*}
In terms of the orthonormal basis $\{u_n\}_{n\in\NN}$ for the subspace
$\mathcal{H}^{1,0}(\PP^1\setminus\Omega)$, the Bergman reproducing kernel can be
written as
\begin{equation*}
K(z,\bar{w})=\sum_{n=1}^{\infty}u_n(z)\otimes \overline{u_n(w)},
\end{equation*}
and does not depend on the choice of the basis. In terms of the Green's function the
Bergman kernel is given by
\begin{equation*}
K(z,\bar{w})=-\frac{\partial^2 G_{DBC}(z,w)}{\partial z\partial\bar{w}} dz\otimes
d\bar{w},
\end{equation*}
and in terms of the conformal map $G$ it has the form
\begin{equation*}
K(z,\bar{w})
=\frac{G^\prime(z)\overline{G^\prime(w)}}{\pi(1-G(z)\overline{G(w)})^2}dz\otimes
d\bar{w} =\sum_{n=1}^\infty u_n(z)\otimes \overline{u_n(w)},
\end{equation*}
where
\begin{equation*}
u_n(z)=\sqrt{\frac{n}{\pi}}G^{-n-1}(z)G^\prime(z)dz,\,\, n\in\NN.
\end{equation*}

\subsection{Deformation theory} Let $G$ be the conformal map $G:\CC\setminus\Omega
\rightarrow\CC\setminus D$, normalized as in the previous section.
Following~\cite{Faber} (see also~\cite{Duren}), we define the analog of Faber
polynomials associated with $G$ by the following Laurent expansion at
$z=\infty~\text{in the region}~|G(z)|>|w|$
\begin{equation*}
\frac{z G^\prime(z)}{G(z)-w}=\sum_{n=0}^\infty F_n(w)z^{-n},
\end{equation*}
obtained by substituting Laurent series for $G(z)$
\begin{equation*}
G(z)=b_{-1}z + b_0 + \frac{b_1}{z} +\dotsb
\end{equation*}
into the geometric series for $(G(z)-w)^{-1}$. We call the degree $n$ polynomials
$F_n$ Faber polynomials of the conformal map $G$. In terms of the inverse map
$g=G^{-1}$ the Faber polynomials can be written explicitly
\begin{equation*}
F_n(w)=[g^n(w)]_{+},
\end{equation*}
where $[g^n]_+$ is a polynomial part of the Laurent series for $g^n$. Faber
polynomials are uniquely characterized by the property
\begin{equation*}
F_n(G(z))=z^n +O(z^{-1})~\text{as}~z\rightarrow\infty.
\end{equation*}
All these facts can be obtained from the representation
\begin{equation*}
F_n(w)=g^n(w)+\frac{1}{2\pi i}\int_{S^1}\frac{g^n(z)}{z-w}dz,\, w\in\CC\setminus D,
\end{equation*}
which follows from the Cauchy integral formula.

\begin{remark} In complex analysis Faber polynomials $P_n$ are usually introduced for
the inverse map $g=G^{-1}:\CC\setminus D\rightarrow\CC\setminus\Omega$ through the
expansion~\cite{Faber,Duren}
\begin{equation*}
\frac{z g^\prime(z)}{g(z)-w}=\sum_{n=0}^\infty P_n(w)z^{-n},
\end{equation*}
where $|g(z)| > |w|$.
\end{remark}

Following~\cite{W-Z,K} (see also references therein) introduce the harmonic moments of
exterior and interior of the contour $C$ by the following formulas
\begin{alignat*}{2}
t_n & =\frac{1}{2\pi in}\int_Cz^{-n}\,\bar{z}\,dz,~n\in\NN, & \qquad t_0 & =
\frac{1}{2\pi i}\int_C\bar{z}\,dz=\frac{1}{\pi} A(\Omega), \\ \intertext{and} v_n & =
\frac{1}{2\pi i}\int_C z^n\,\bar{z}\,dz,~n\in\NN, & \qquad v_0 & = \frac{2}{\pi}
\int_{\Omega}\log|z|\,\bar{z}\,d^2z.
\end{alignat*}
Here $A(\Omega)$ is the area of the interior domain $\Omega$ with respect to the
standard Lebesgue measure $d^2z=|dz\wedge d\bar{z}|/2$ on $\CC$. According
to~\cite{K,W-Z} and references therein, parameters $\{t_0,t_n,\bar{t}_n\}_{n\in\NN}$
are local coordinates for the space $\mathcal{C}$ in some neighborhood of the contour
$C$. Consequently, there is a foliation of $\mathcal{C}$ with the leaves
$\Tilde{\mathcal{C}}_a$ of contours of fixed area $a>0$. The leaves
$\Tilde{\mathcal{C}}_a$ can be considered (at least locally) as infinite-dimensional
complex manifolds with complex coordinates $\{t_n\}_{n\in\NN}$.

\begin{remark} There is also a foliation of $\mathcal{C}$ with the leaves
$\mathcal{C}_r$ of contours of conformal radius $r>0$. For fixed $r$ the coordinate
$t_0$ can be expressed in terms of $t_n,\bar{t}_n$ provided that $\partial r/\partial
t_0\neq 0$ so that the leaves $\mathcal{C}_r$ can be considered (at least locally) as
infinite-dimensional complex manifolds. Similarly, there is a foliation of
$\mathcal{C}$ with the leaves $\mathcal{C}^d$ of contours of transfinite diameter
$d>0$. In this paper we do not address a very interesting problem of global
description of complex manifolds $\mathcal{C}_r~\text{and}~\mathcal{C}^d$.
\end{remark}

The goal of the deformation theory is to describe the tangent vector space
$T_C\mathcal{C}$ to the manifold $\mathcal{C}$ at a contour $C$ in terms of the data
associated with $C$. By definition, a deformation of the contour $C$ is a smooth
family of contours $\{C_t\}_{t\in(-\epsilon,\epsilon)}$ such that $C_0=C$. The
smoothness property is a condition that there exists a parameterization of contours
$C_t$ which depends smoothly on $t$: $C_t=\{z(\sigma,t),\, \sigma\in\RR/2\pi\ZZ\}$ for
every $|t|<\epsilon$, where $z(\sigma,t)\in
C^\infty(\RR/2\pi\ZZ\times(-\epsilon,\epsilon))$. Corresponding infinitesimal
deformation is defined as the vector field $v=\dot{z}(\sigma)d/d\sigma$ along $C$,
where dot stands for $d/dt|_{t=0}$,
\begin{equation*}
\dot{z}(\sigma):=\frac{\partial z(\sigma,t)}{\partial t}\bigg|_{t=0}.
\end{equation*}

By definition, a deformation $C_t$ is called trivial if it consists of
reparameterizations of the contour $C$, so that the vector field $v$ is tangential to
$C$. The deformation $C_t$ is called infinitesimally trivial if the corresponding
vector field $v$ is tangential to $C$. The tangent vector space $T_C\mathcal{C}$ can
be canonically identified with the quotient space of the real vector space of all
vector fields along $C$ by the subspace of vector fields tangential to $C$.
Equivalently, $T_C\mathcal{C}$ is a real vector space of normal vector fields to $C$.

With every deformation $C_t$ one can associate the following 1-form on $C$
\begin{equation*}
\dot{\omega}_C:=\frac{\partial\bar{z}}{\partial t}\bigg|_{t=0}dz - \frac{\partial
z}{\partial t}\bigg|_{t=0}d\bar{z}.
\end{equation*}
The 1-form $\dot{\omega}_C$ is a restriction to the contour $C$ of a $d^{-1}$ of the
Lie derivative of the standard 2-form $d\bar{z}\wedge d z$ on $\Omega$. It satisfies
the following ``calculus formula''
\begin{equation*}
\frac{d}{dt}\bigg|_{t=0}\int_{C_t}f(z,\bar{z},t)dz=\int_C(\dot{f} dz+ \frac{\partial
f}{\partial\bar{z}}\dot{\omega}_C)
\end{equation*}
for any smooth function $f$ in the domain $(\bigcup_{t\in (-\epsilon,\epsilon)}C_t)
\times (-\epsilon,\epsilon)$.
\begin{remark} Classically, as it goes back to Volterra and Hadamard,
deformations of a contour $C$ are described by
\begin{equation*}
z(\sigma,t)=z(\sigma)+t\delta n(\sigma)n(\sigma),~\sigma\in\RR/2\pi\ZZ,
\end{equation*}
where $z(\sigma)$ defines the contour $C$, $\delta n(\sigma)\in
C^\infty(\RR/2\pi\ZZ,\RR)$, and $n(\sigma)$ is the outer normal to $C$. For $t$
sufficiently small this equation defines smooth closed Jordan curve $C_t$. The
relation between two approaches is given by the elementary calculus
\begin{equation*}
\delta n ds=\frac{1}{2i}\dot{\omega}_C,
\end{equation*}
where $ds:=|z^\prime(\sigma)|d\sigma$ is a 1-form on $C$.
\end{remark}

The basic facts of the deformation theory of contours are summarized in the following
theorem. Parts (i) and (ii) were proved in~\cite{K,W-Z} and we refer to them as
``Krichever's lemma''.

\begin{theorem}

\begin{itemize}

\item[(i)] Any deformation $C_t$ of the contour $C$
which does not change harmonic moments of exterior $\{t_0,t_n,\bar{t}_n\}$ is
infinitesimally trivial. The parameters
$\{t_0-t_0(C),t_n-t_n(C),\bar{t}_n-\bar{t}_n(C)\}$ are local coordinates on
$\mathcal{C}$ in some neighborhood of $C$.
\item[(ii)] The following $1$-forms on $C$
\begin{equation*}
\dot{\omega}^{(n)}_C:=\frac{\partial\bar{z}}{\partial t_n}\bigg|_{t_n=t_n(C)}dz -
\frac{\partial z}{\partial t_n}\bigg|_{t_n=t_n(C)}d\bar{z},\,n\in\{0\}\cup\NN,
\end{equation*}
extends to meromorphic $(1,0)$-forms on the double $\PP^1_C$ of the exterior domain
$\PP^1\setminus\Omega$ with a single pole at $\infty$ of order $n+1$ if $n\in\NN$ and
simple poles at $\infty$ and $\overline{\infty}$ with residues $1$ and $-1$ if $n=0$.
Explicitly
\begin{equation*}
\dot{\omega}^{(n)}_C = d(F_n\circ G),~\dot{\omega}^{(0)}_C = d\log G
\end{equation*}
in the domain $\PP^1\setminus\Omega$ and
\begin{equation*}
\dot{\omega}^{(n)}_C = d(F_n\circ 1/\overline{G}),~\dot{\omega}^{(0)}_C  =
d\log1/\overline{G}
\end{equation*}
in the domain $\overline{\PP^1\setminus\Omega}$, where $F_n$ is the $n$-th Faber
polynomial of $G$, $n\in\NN$.

\item[(iii)] The $1$-forms $\dot{\omega}^{(n)}_C$ satisfy the property
\begin{equation*}
\frac{1}{2\pi i}\int_C\frac{z^{-m}}{m}\,\dot{\omega}^{(n)}_C=\delta_{mn},
\end{equation*}
and can be identified with the vector fields $\partial/\partial t_n$. For every $a>0$
the holomorphic tangent vector space $T^\prime_C\Tilde{\mathcal{C}}_a$ to
$\Tilde{\mathcal{C}}_a$ at $C$ is canonically identified with the complex vector space
$\mathcal{M}^{1,0}(\PP^1_C)$ of meromorphic $(1,0)$-forms on $\PP^1_C$ with a single
pole at $\infty$ of order $\geq 2$.
\item[(iv)] The holomorphic cotangent vector space $T^{\prime^\ast}_C\Tilde{\mathcal{C}}_a$
to $\Tilde{\mathcal{C}}_a $ at $C$ is naturally identified with the complex vector
space $\mathcal{H}^{1,0}(\PP^1\setminus\Omega)$ of all holomorphic $(1,0)$-forms on
$\PP^1\setminus\Omega$ which are smooth up to the boundary, with the pairing
\begin{equation*}
(~,~)_C: T^\prime_C\Tilde{\mathcal{C}}_a\otimes
T^{\prime^\ast}_C\Tilde{\mathcal{C}}_a\rightarrow\CC
\end{equation*}
given by
\begin{equation*}
(\omega,u)_C:= \frac{1}{2\pi i}\int_C d^{-1}u\,\omega.
\end{equation*}
Differentials $dt_n$ correspond to $(1,0)$-forms
\begin{equation*}
dt_n(z):=d(z^{-n}/n)=-z^{-n-1}dz.
\end{equation*}
\end{itemize}
\end{theorem}
\begin{proof} As in~\cite{K,W-Z}, we start with the following Riemann-Hilbert
problem. Find functions $S_{+}$ and $S_-$ that are holomorphic in the domains $\Omega$
and $\CC\setminus\Omega$ respectfully, are smooth in the corresponding closed domains
and on their common boundary $C$ satisfy
\begin{equation*}
S_{+}(z)-S_{-}(z)=\bar{z}~\text{for all}~z\in C.
\end{equation*}
With the normalization $S_-(\infty)=0$ this problem has a unique solution given by the
Cauchy integral
\begin{equation*}
S_{\pm}(z)=\frac{1}{2\pi i}\int_C\frac{\bar{w}dw}{w-z},
\end{equation*}
where $z\in\Omega$ for the $+$ sign and $z\in\CC\setminus\Omega$ for the $-$ sign.

Now let $C_t$ be a deformation of $C$ satisfying conditions in part (i). By the
calculus formula,
\begin{equation*}
\dot{S}_{\pm}(z)=\frac{1}{2\pi i}\int_C\frac{\dot{\bar{w}} dw-\dot{w}d\bar{w}}{w-z},
\end{equation*}
so that on $C$
\begin{equation*}
\dot{\omega}_C=\dot{S}_{+}dz-\dot{S}_{-}dz.
\end{equation*}
Since in some neighborhood of $0$
\begin{equation*}
S_+(z)=\sum_{n=1}^\infty n t_n z^{n-1},
\end{equation*}
we see that if $dt_n/dt|_{t=0}=0$ for all $n\in\NN$ then $\dot{S}_+=0$ in $\Omega$.
This implies that on $C$
\begin{equation*}
\dot{\omega}_C=-\dot{S}_{-}dz|_C,
\end{equation*}
and $\dot{\omega}_C$ admits holomorphic continuation as a $(1,0)$-form on
$\CC\setminus\Omega$. Since $dt_0/dt|_{t=0}=0$ and
\begin{equation*}
S_-(z)=-\frac{t_0}{z} + O(z^{-2})~\text{as}~z\rightarrow\infty,
\end{equation*}
we get $\dot{S}_{-}(z)=O(z^{-2})$, so that $\dot{\omega}_C$ is also regular at
$\infty$. Since, by definition, the 1-form $\dot{\omega}_C$ is pure imaginary on $C$,
by Riemann-Schwarz reflection principle it can be analytically continued to a
holomorphic $(1,0)$-form on $\PP^1_C$. Since $\PP^1_C$ has genus 0 we conclude that
$\dot{\omega}_C=0$. In particular, $\dot{\omega}_C=0$ on $C$, which by Remark 1.4 is
equivalent to the condition that vector field $v$ corresponding to the deformation
$C_t$ is tangential to $C$. This proves the part (i).

For the proof of the part (ii), set $t=\re(t_n - t_n(C))$. We have for $n\in\NN$
\begin{equation*}
\dot{S}_+=nz^{n-1},
\end{equation*}
so that the $1$-form $\dot{\omega}_C-dz^n$ admits holomorphic continuation to the
domain $\CC\setminus\Omega$ and is regular at $\infty$. Similarly, for $n=0$ the
$1$-form $\dot{\omega}_C$ is holomorphic on $\CC\setminus\Omega$ with a simple pole at
$\infty$ with residue $1$. As before, we conclude that 1-form $\dot{\omega}_C$ admits
a meromorphic continuation to $\PP^1_C$ with only poles at $\infty$ and
$\overline{\infty}$. For $n=0$ they are simple poles with residues $1$ and $-1$, so
that using the global coordinate $\zeta$ on $\PP_C^1$ we get
\begin{equation*}
\dot{\omega}^{(0)}_C=\frac{d\zeta}{\zeta}.
\end{equation*}
In particular, $\dot{\omega}^{(0)}_C=d\log G$ on $\PP^1\setminus\Omega$ and
$\dot{\omega}^{(0)}_C=d\log 1/\overline{G}$ on $\overline{\PP^1\setminus\Omega}$. For
$n\in\NN$ we get
\begin{equation*}
\dot{\omega}_C=d(z^n + O(z^{-1}))~\text{as}~z\rightarrow\infty,~\text{and}~
\dot{\omega}_C=-d(\bar{z}^n + O(z^{-1}))~\text{as}~z\rightarrow\bar\infty.
\end{equation*}

Similarly, setting $t=\im(t_n -t_n(C))$ we get for $n\in\NN$
\begin{equation*}
\dot{\omega}_C=id(z^n + O(z^{-1}))~\text{as}~z\rightarrow\infty~\text{and}~
\dot{\omega}_C=id(\bar{z}^n + O(z^{-1}))~\text{as}~z\rightarrow\bar\infty.
\end{equation*}
From here we conclude that for $n\in\NN$ the $1$-forms $\dot{\omega}^{(n)}_C$ are
meromorphic on $\PP^1_C$ with the only pole at $\infty$ and
\begin{equation*}
\dot{\omega}^{(n)}_C=d(z^n + O(z^{-1}))~\text{as}~z\rightarrow\infty.
\end{equation*}

Similarly, the 1-forms
\begin{equation*}
\dot{\omega}^{(\bar{n})}_C:=\frac{\partial\bar{z}}{\partial \bar{t}_n}\bigg|_{t=0}dz -
\frac{\partial z}{\partial \bar{t}_n}\bigg|_{t=0}d\bar{z}
\end{equation*}
are meromorphic on $\PP^1_C$ with the only pole at $\overline{\infty}$ and
\begin{equation*}
\dot{\omega}^{(\bar{n})}_C=-d(\bar{z}^n +
O(\bar{z}^{-1}))~\text{as}~z\rightarrow\bar\infty.
\end{equation*}
Using characteristic property of the Faber polynomials $F_n$ we get that on
$\PP^1\setminus\Omega$
\begin{equation*}
\dot{\omega}^{(n)}_C=d(F_n\circ G)~\text{and}~
\dot{\omega}^{(\bar{n})}_C=-d(\overline{F_n\circ 1/\overline{G}}),
\end{equation*}
and on $\overline{\PP^1\setminus\Omega}$
\begin{equation*}
\dot{\omega}^{(n)}_C=d(F_n\circ 1/\overline{G})~\text{and}~
\dot{\omega}^{(\bar{n})}_C=-d(\overline{F_n\circ G}).
\end{equation*}

The proof of parts (iii)-(iv) is now straightforward. By definition of the harmonic
moments $t_n$, the calculus formula and the Cauchy theorem we get
\begin{equation*}
\frac{\partial t_n}{\partial t_m}=\frac{1}{2\pi i n}\int_C z^{-n}\dot{\omega}^{(m)}_C
=\delta_{mn},
\end{equation*}
so that 1-forms $\dot{\omega}^{(m)}_C$ correspond to the vector fields
$\partial/\partial t_m$ on $\mathcal{C}$. In particular, the same formula shows that
corresponding differentials $dt_n$ --- $(1,0)$-forms on $\mathcal{C}$, are given by
the 1-forms $d(z^{-n}/n)$ on $C$ and satisfy
\begin{equation*}
\left(\frac{\partial}{\partial t_m}\,, dt_n\right)_C=\delta_{mn}.
\end{equation*}
\end{proof}

\begin{remark} Let $\mathbf{d}=\mathbf{d}^\prime + \mathbf{d}^{\prime\prime}$
be the decomposition of the de Rham differential on $\Tilde{\mathcal{C}}_a$ into
$(1,0)$ and $(0,1)$ components with respect to the complex structure defined by the
harmonic moments $t_n$. According to part (iv) of the theorem, for any smooth function
$F$ on $\Tilde{\mathcal{C}}_a$ the $(1,0)$-form
\begin{equation*}
\mathbf{d}^\prime F:=\sum_{n=1}^{\infty}\frac{\partial F}{\partial t_n} dt_n
\end{equation*}
at $C\in\Tilde{\mathcal{C}}_a$ can be identified with the holomorphic $(1,0)$-form on
$\PP^1\setminus\Omega$ defined by the following Laurent expansion at $z=\infty$
\begin{equation*}
-\sum_{n=1}^{\infty}z^{-n-1}\frac{\partial F}{\partial t_n} dz.
\end{equation*}
The convergence of this series and holomorphy of the corresponding $(1,0)$-form on
$\PP^1\setminus\Omega$ follow from the smoothness of the function $F$ --- the
existence of $\mathbf{d}F$. An example of smooth function $F$ is given by
\begin{equation*}
F=\frac{1}{2\pi i}\int_C h(z)\bar{z} dz,
\end{equation*}
where $h$ is holomorphic on $\PP^1\setminus\Omega$ and is smooth up the boundary. In
this case,
\begin{equation*}
\mathbf{d}^\prime F=h^\prime(z) dz.
\end{equation*}
Though we do not address here the question of defining various functional classes on
$\Tilde{\mathcal{C}}_a$, all said above holds for the class of real-analytic
functions.
\end{remark}

Here we introduce a natural Hermitian metric on complex manifolds
$\Tilde{\mathcal{C}}_a$, which turns out to be K\"{a}hler, as we shall prove in
Section 3. Namely, for every $C\in\Tilde{\mathcal{C}}_a$ consider the following inner
product in the holomorphic tangent vector space $T^{\prime}_C\Tilde{\mathcal{C}}_a$
\begin{equation*}
H\left(\frac{\partial}{\partial t_m},\frac{\partial}{\partial
t_n}\right)=h^{m\bar{n}}:= - \frac{1}{(2\pi i)^2}\int_{C_+}\int_{C_+} z^m \bar{w}^n
K(z,\bar{w}).
\end{equation*}
where $C_+$ is an arbitrary contour containing $C$ inside (not that $K$ is singular as
$z=w\in C$). From the representation of the Bergman kernel in terms of the orthonormal
system in Section 1.2 it easily follows that this inner product is positive-definite.
Using this orthonormal system it is not difficult to show that components
$h^{m\bar{n}}$ are smooth on $\Tilde{\mathcal{C}}_a$ and thus define the Hermitian
metric $H$. In terms of this metric we get the following Laurent expansion at
$z=w=\infty$ for the Bergman reproducing kernel $K$,
\begin{equation*}
K(z,\bar{w})=\sum_{m,n=1}^\infty h^{m\bar{n}} dt_m(z)\otimes\overline{dt_n(w)}.
\end{equation*}
\begin{remark}
It is instructive to compare the complex structure on $\mathcal{C}_1$ introduced by
A.A.~Kirillov~\cite{Kirillov1} with the complex structure on $\Tilde{\mathcal{C}}_1$
defined (at least locally) by the harmonic moments of exterior. According
to~\cite{N-V} (see Remark 1.1) the former is a pull-back by the inclusion map of the
Ahlfors-Bers complex structure on $T(1)$ and is defined using quadratic (or,
equivalently, Beltrami) differentials, whereas the latter is defined using
$(1,0)$-forms --- the ordinary differentials. Correspondingly, the Hermitian metric on
$\Diff_+(S^1)/\Mob(S^1)$ is given by the Petersson inner product of holomorphic
quadratic (equivalently, Beltrami) differentials~\cite{Bers,N-V}, whereas the
Hermitian metric on $\Tilde{\mathcal{C}}_1$ is defined via the canonical inner product
of holomorphic $(1,0)$-forms on $\PP^1\setminus\Omega$.
\end{remark}

\section{Bosonic action functional and partition function}

For a classical field $\varphi\in C^\infty(\PP^1,\RR)$ consider the following action
functional
\begin{equation*}
S_0(\varphi):=\frac{1}{8}\int_{\PP^1}d\varphi\wedge \ast d \varphi=
\frac{i}{4}\int_{\PP^1}\partial\varphi\wedge\bar{\partial}\varphi,
\end{equation*}
which describes the standard theory of free bosons on the Riemann sphere $\PP^1$.
Corresponding partition function is defined by the functional integral
\begin{equation*}
\langle\mathbf{1}\rangle_0:=\int_{C^\infty(\PP^1,\RR)/\RR}[\mathcal{D}\varphi]
e^{-\frac{1}{\pi}S_0(\varphi)},
\end{equation*}
where integration goes over the coset $C^\infty(\PP^1,\RR)/\RR$ and reflects the
symmetry $\varphi\mapsto\varphi+c$. As in Part I, mathematically rigorous definition
requires a choice of a conformal metric $ds^2$ on $\PP^1$ and leads to the result
\begin{equation*}
\langle\mathbf{1}\rangle_0=\int_{C^\infty(\PP^1,\RR)/\RR}[\mathcal{D}\varphi]
e^{-\frac{1}{2\pi}\int_{\PP^1}\Delta_0\varphi\ast\varphi}=
\left(\frac{\Area(\PP^1)}{\det_\zeta\Delta_0}\right)^{1/2}.
\end{equation*}
Here the area term is a contribution from zero modes --- the one-dimensional kernel of
the $\bar\partial$-Laplacian $\Delta_0$ of the metric $ds^2$ acting on functions on
$\PP^1$.

For every $C\in\mathcal{C}$ introduce the following analog of the topological term
(cf.~with the discussion in Section 2 of Part I):
\begin{equation*}
S_{top}(\varphi):=\int_{\CC}(A(\Omega)\delta_0-\chi_\Omega)\varphi d^2z=
A(\Omega)\varphi(0)-\int_{\Omega}\varphi d^2z.
\end{equation*}
Here $\chi_\Omega$ is a characteristic function of the domain $\Omega$, and $\delta_0$
is a Dirac delta-function at $0$ with respect to the Lebesgue measure $d^2z$. The
functional $S_{top}$ has the property $S_{top}(\varphi+c)=S_{top}(\varphi)$.

The total bosonic action
\begin{equation*}
S_C(\varphi):=S_0(\varphi)+S_{top}(\varphi)
\end{equation*}
defines the theory of free bosons on $\PP^1$ in the presence of a contour $C$ and we
consider a family of such field theories parameterized by $\mathcal{C}$.

For every $C\in\mathcal{C}$ define the partition function of the corresponding quantum
field theory by the following functional integral
\begin{equation*}
\langle\mathbf{1}\rangle_C:=\int_{C^\infty(\PP^1,\RR)/\RR}
[\mathcal{D}\varphi]e^{-\frac{1}{\pi}S_C(\varphi)}.
\end{equation*}

Mathematically rigorous definition is the following. Approximate, in the
distributional sense, characteristic function $\chi_\Omega$ and Dirac delta-function
$\delta_0$ by smooth functions $\chi^{(\epsilon)}_\Omega$ and $\delta^{(\epsilon)}_0$
with compact supports satisfying
\begin{equation*}
\int_{\CC} (A(\Omega)\delta^{(\epsilon)}_0-\chi^{(\epsilon)}_\Omega)d^2z=0
\end{equation*}
and define
\begin{align*}
\langle\mathbf{1}\rangle_C:= & \lim_{\epsilon\rightarrow 0}\exp\{\frac{A^2(\Omega)}
{\pi^2}\int_\CC\int_\CC\log|z-w|\delta^{(\epsilon)}_0(z)\delta^{(\epsilon)}_0(w)
d^2zd^2w\} \\ &
\int_{C^\infty(\PP^1,\RR)/\RR}[\mathcal{D}\varphi]e^{-\frac{1}{\pi}S^{(\epsilon)}_C(\varphi)},
\end{align*}
where
\begin{equation*}
S^{(\epsilon)}_C(\varphi):=S_0(\varphi)+\int_{\CC}
(A(\Omega)\delta^{(\epsilon)}_0-\chi^{(\epsilon)}_\Omega)\varphi d^2z.
\end{equation*}

We introduce the $\tau$-function $\tau=\tau(C)$ of the smooth Jordan contour $C$ as
the \emph{normalized expectation value} of $C$, defined as follows
\begin{equation*}
\tau=\langle\langle C\rangle\rangle:=
\frac{\langle\mathbf{1}\rangle_C}{\langle\mathbf{1}\rangle_0}.
\end{equation*}
\begin{proposition}
The $\tau$-function of the contour $C$ is well-defined and is given explicitly by the
following expression
\begin{align*}
\log\tau & =-\frac{1}{\pi^2}\int_\Omega\int_\Omega \log|z-w|d^2z d^2w+\frac{2}{\pi^2}
A(\Omega)\int_\Omega \log|z|d^2z \\ & =-\frac{1}{\pi^2}\int_\Omega\int_\Omega
\log|\frac{1}{z}-\frac{1}{w}|d^2z d^2w.
\end{align*}
\end{proposition}
\begin{proof} It is another standard computation. Consider the Gaussian integral
\begin{equation*}
\langle\mathbf{1}\rangle_C^{(\epsilon)}:=\int_{C^\infty(\PP^1,\RR)/\RR}
[\mathcal{D}\varphi]e^{-\frac{1}{\pi}S^{(\epsilon)}_C(\varphi)}.
\end{equation*}
and make the change of variables $\varphi=\Phi^{(\epsilon)}+\tilde\varphi$, where
$\Phi^{(\epsilon)}$ is uniquely determined by the condition that
$S^{(\epsilon)}_C(\Phi^{(\epsilon)}+\tilde\varphi)$ does not contain linear terms in
$\tilde\varphi$ and by the normalization $\Phi^{(\epsilon)}(\infty)=0$. Using the
Stokes theorem we get
\begin{equation*}
-\frac{\partial^2\Phi^{(\epsilon)}(z)}{\partial
z\partial\bar{z}}=\lambda^{(\epsilon)}(z),
\end{equation*}
where $\lambda^{(\epsilon)}:=\chi^{(\epsilon)}_\Omega-A(\Omega)\delta^{(\epsilon)}_0$.
The function $\lambda^{(\epsilon)}$ is smooth, has compact support and
\begin{equation*}
\int_{\CC}\lambda^{(\epsilon)}(z)d^2z=0,
\end{equation*}
so that
\begin{equation*}
\Phi^{(\epsilon)}(z)=-\frac{2}{\pi}\int_\CC\log|z-w|\lambda^{(\epsilon)}(w)d^2w.
\end{equation*}
Since
\begin{equation*}
S^{(\epsilon)}_C(\Phi^{(\epsilon)}+\tilde\varphi)=S_0(\tilde\varphi)+\frac{1}{2}\int_{\CC}
\Phi^{(\epsilon)}(z)\lambda^{(\epsilon)}(z)d^2z
\end{equation*}
we finally obtain
\begin{equation*}
\langle\mathbf{1}\rangle_C^{(\epsilon)}=\langle
\mathbf{1}\rangle_0\exp\{-\frac{1}{\pi^2}\int_\CC\int_\CC
\log|z-w|\lambda^{(\epsilon)}(z)\lambda^{(\epsilon)}(w)d^2zd^2w\}.
\end{equation*}
Multiplying by the regularization factor and passing to the limit $\epsilon\rightarrow
0$ we see that $\log\tau_C$ is well-defined and is given by the formula above.
\end{proof}
\begin{corollary} The $\tau$-function is $-1/\pi^2$ times regularized energy of the
pseudo-measure $d\mu=d^2z-A(\Omega)\delta_0$ on the domain $\Omega$, where $d^2z$ is
the Lebesgue measure and $\delta_0$ --- the delta-measure at 0.
\end{corollary}
\begin{proof}
Indeed, the energy $I(\nu)$ of a Borel measure $d\nu$ on $\Omega$ is defined by (see,
e.g.~\cite{potential})
\begin{equation*}
I(\nu):=\int_\Omega\int_\Omega \log|z-w|d\nu(z) d\nu(w).
\end{equation*}
In our case, due to the presence of a delta-measure, we formally have
$I(\mu)=-\infty$. However, with the above regularization $I(\mu)=-\pi^2 \tau$.
\end{proof}

\begin{remark}
It follows from the Proposition 2.1 that the $\tau$-function $\tau_C$ coincides with
Mineev-Weinstein--Wiegmann--Zabrodin $\tau$-function $\tau_{WZ}$ for the analytic
closed Jordan curve
$C$, introduced in~\cite{W-Z} and computed in~\cite{K-K-MW-W-Z}! It is interesting to
compare these two approaches. Specifically, in~\cite{W-Z, K-K-MW-W-Z} the
$\tau$-function $\tau_{WZ}$ appears as a dispersionless limit of the Hirota's
$\tau$-function for the integrable two-dimensional Toda hierarchy and also as a large
$N$ limit of a partition function of a certain ensemble of random $N\times N$
matrices. In our approach the $\tau$-function $\tau_C$ is just a partition function of
a quantum field theory of free bosons on $\PP^1$ parameterized by a smooth contour
$C\in\mathcal{C}$. This is quite analogous to the definition of the $\btau$-function
in Part I as a partition function of a quantum field theory of free bosons on $X$
parameterized by a holomorphic line bundle $L\in\Pic^{g-1}(X)\setminus\mathcal{E}$.
\end{remark}

We also consider the quantum theory of free bosons in the exterior domain
$\PP^1\setminus\Omega$, defined by the following functional
\begin{equation*}
S_{ext}(\varphi):=\frac{1}{8}\int_{\PP^1\setminus\Omega}d\varphi\wedge \ast d \varphi=
\frac{i}{4}\int_{\PP^1\setminus\Omega}\partial\varphi\wedge\bar{\partial}\varphi,
\end{equation*}
where $\varphi\in C^\infty_{DBC}(\PP^1\setminus\Omega)$
and satisfies the Dirichlet boundary condition
\begin{equation*}
\varphi|_{C}=0.
\end{equation*}
The corresponding partition function is defined by the functional integral
\begin{equation*}
\langle\mathbf{1}\rangle_{DBC}:=\int_{C^{\infty}_{DBC}(\PP^1\setminus\Omega,\RR)}
[\mathcal{D}\varphi]e^{-\frac{1}{\pi}S_{ext}(\varphi)}
\end{equation*}
and is given explicitly as follows
\begin{equation*}
\langle\mathbf{1}\rangle_{DBC}=\left(\frac{\Area(\PP^1\setminus\Omega)}{\det_\zeta\Delta_0}
\right)^{1/2},
\end{equation*}
where $\Delta_0$ is a $\bar\partial$-Laplacian of the metric $ds^2$ acting on
functions on $\PP^1\setminus\Omega$ satisfying the Dirichlet boundary condition.

\section{Current Ward identities}
Introduce holomorphic and anti-holomorphic components of the bosonic field current
$d\varphi$ as follows
\begin{equation*}
d\varphi=\jmath + \bar{\jmath},
\end{equation*}
where $\jmath=\partial\varphi~\text{and}~\bar{\jmath}=\bar\partial\varphi$. By
definition, multi-point correlation functions of current components are given by the
following functional integral
\begin{multline*}
\langle \jmath(z_1)\dotsb \jmath(z_m)\bar{\jmath}(w_1)\dotsb\bar{\jmath}(w_n)\rangle
\\ :=\int_{C^\infty(\PP^1,\RR)/\RR}[\mathcal{D}\varphi] \,\jmath(z_1)\dotsb
\jmath(z_m)\bar{\jmath}(w_1)\dotsb\bar{\jmath}(w_n)\,e^{-\frac{1}{\pi}S_C(\varphi)}.
\end{multline*}
Correlation functions for the standard theory of free bosons on $\PP^1$ with action
functional $S_0$ (which formally corresponds to $C=\emptyset$ -- the empty set) and
for the theory on $\PP^1\setminus\Omega$ with action functional $S_{ext}$ and
Dirichlet boundary condition are defined similarly. We denoted them by
$\langle\dotsb\rangle_0$ and $\langle\dotsb\rangle_{DBC}$ respectfully.

We start with the computation of the normalized 1-point correlation function of the
holomorphic component $\jmath$, defined by
\begin{equation*}
\langle\langle \jmath(z)\rangle\rangle:=\frac{\langle
\jmath(z)\rangle}{\langle\mathbf{1}\rangle_C}.
\end{equation*}
Repeating the proof of the Proposition 2.1, which is based on the change of variables
$\varphi=\Phi^{(\epsilon)}+ \tilde\varphi$, and using the standard fact that $\langle
\tilde\jmath(z) \rangle_0=0$, where $\tilde\jmath(z)=\partial\tilde\varphi(z)$
(cf.~Section 4 in Part I) we get
\begin{equation*}
\int_{C^\infty(\PP^1,\RR)/\RR}[\mathcal{D}\varphi]
\,\jmath(z)\,e^{-\frac{1}{\pi}S^{(\epsilon)}_C(\varphi)}= \langle\mathbf{1}
\rangle_C^{(\epsilon)}\frac{\partial\Phi^{(\epsilon)}(z)}{\partial z}dz.
\end{equation*}
Passing to the limit $\epsilon\rightarrow 0$ we obtain
\begin{equation*}
\langle\langle \jmath(z)\rangle\rangle= \frac{\partial \Phi(z)}{\partial z}dz,
\end{equation*}
where $\Phi:=\lim_{\epsilon\rightarrow 0}\Phi^{(\epsilon)}$ and is given explicitly by
\begin{equation*}
\Phi(z)=\frac{2 A(\Omega)}{\pi}\log|z| - \frac{2}{\pi}\int_\Omega\log|z-w|d^2w.
\end{equation*}
The function $\Phi(z)$ can be characterized as a continuous solution of the equation
\begin{equation*}
-\frac{\partial^2\Phi(z)}{\partial z\partial\bar{z}}=
\begin{cases}
\chi_\Omega(z) - A(\Omega)\delta_0(z)& \text{if $z\in\Omega$,} \\ 0 & \text{if
$z\in\CC\setminus\Omega$,}
\end{cases}
\end{equation*}
normalized by the condition $\Phi(\infty)=0$, and is a logarithmic potential of the
pseudo-measure $d\mu=d^2z - A(\Omega)\delta_0$ on $\Omega$. It also follows from the
integral representation that holomorphic on $\PP^1\setminus\Omega$ function
$\partial\Phi/\partial z$ coincides with the function $S_-(z)+t_0/z$ (see Section
1.3), and has the following Laurent expansion at $z=\infty$
\begin{equation*}
\frac{\partial \Phi(z)}{\partial z}=-\sum_{n=1}^\infty v_n z^{-n-1}.
\end{equation*}

On the other hand, consider the partial derivative
\begin{equation*}
\frac{\partial \log\tau}{\partial t_n}=\frac{1}{\langle\mathbf{1}\rangle_C}
\frac{\partial\langle\mathbf{1}\rangle_C}{\partial t_n} ,\,\,n\in\NN.
\end{equation*}
We can evaluate it by differentiating under the functional integral sign, which can be
easily justified using the rigorous definition of the partition function
$\langle\mathbf{1}\rangle_C$. Since only topological term in the action functional
depends on the domain $\Omega$, the computation is based on another calculus formula
\begin{equation*}
\frac{\partial}{\partial t_n}\int_{\Omega}\varphi d^2z=\frac{1}{2i}\int_C
\varphi\dot{\omega}^{(n)}_C.
\end{equation*}
Arguing as in the proof of Proposition 2.1 we get
\begin{equation*}
\frac{\partial\log\tau}{\partial t_n} =\frac{1}{2\pi i}\langle\int_C\varphi
\dot{\omega}^{(n)}_C\rangle/\langle\mathbf{1}\rangle_C= \frac{1}{2\pi i}\int_C\Phi
\dot{\omega}^{(n)}_C.
\end{equation*}
Using Krichever's lemma and integration by parts we obtain
\begin{equation*}
\frac{\partial\log\tau}{\partial t_n}=\frac{1}{2\pi i}\int_C\Phi d(F_n\circ G)
=-\frac{1}{2\pi i}\int_C F_n\circ G (\frac{\partial\Phi}{\partial z} dz +
\frac{\partial\Phi}{\partial\bar{z}}d\bar{z}),
\end{equation*}
where $\partial\Phi/\partial z$ and $\partial\Phi/\partial\bar{z}$ are boundary values
on $C$ of holomorphic and anti-holomorphic functions on the exterior domain
$\PP^1\setminus\Omega$. The latter can be easily justified by considering smooth
functions $\Phi^{(\epsilon)}$ first and then passing to the limit $\epsilon\rightarrow
0$. Using the characteristic property of Faber polynomials and Cauchy theorem we have
\begin{equation*}
-\frac{1}{2\pi i}\int_C (F_n\circ G)\frac{\partial\Phi}{\partial z} dz =-\frac{1}{2\pi
i}\int_C \frac{\partial\Phi}{\partial z} z^n dz =v_n.
\end{equation*}
We claim that the second integral in the formula for $\partial\log\tau/\partial t_n$
is 0. Using that $|G|=1$ on $C$, we have
\begin{equation*}
\overline{\int_C (F_n\circ G)\frac{\partial\Phi}{\partial\bar{z}}d\bar{z}} =\int_C
(\bar{F}_n \circ 1/G) \frac{\partial\Phi}{\partial z} dz,
\end{equation*}
where $\bar{F}_n(z):=\overline{F_n(\bar{z})}$. Since $\partial\Phi/\partial
z=O(z^{-2})$ as $z\rightarrow\infty$, the integral indeed vanishes.

Thus we proved the Ward identity for the 1-point correlation function
\begin{equation*}
\frac{\partial\log\tau}{\partial t_n}=-\frac{1}{2\pi i}\int_C z^n
\langle\langle\jmath(z)\rangle\rangle,\,\,n\in\NN.
\end{equation*}

The case $n=0$ can be considered similarly with the only difference that since
$A(\Omega)=\pi t_0$ one needs to differentiate the regularization factor as well. We
have the following computation
\begin{align*}
\frac{\partial\log\tau}{\partial t_0}& =\lim_{\epsilon\rightarrow 0}\biggl(
\frac{1}{2\pi i}\int_C\Phi^{(\epsilon)}\frac{dG}{G}-\int_{\CC}\Phi^{(\epsilon)}(z)
\delta^{(\epsilon)}(z)d^2z \\ &
+\frac{2A(\Omega)}{\pi}\int_\CC\int_\CC\log|z-w|\delta^{(\epsilon)}_0(z)
\delta^{(\epsilon)}_0(w)d^2zd^2w \biggr) \\ &=\frac{1}{2\pi i}\int_C\Phi\frac{dG}{G} +
\lim_{\epsilon\rightarrow 0} \int_{\CC}\biggl(\frac{2A(\Omega)}{\pi}\log|z| -
\Phi^{(\epsilon)}(z)\biggr)\delta^{(\epsilon)}_0(z) d^2z \\ &=\lim_{z\rightarrow
0}\biggl(\frac{2A(\Omega)}{\pi}\log |z| - \Phi(z)\biggr) \\ &=v_0,
\end{align*}
with an easy justification of all steps. Here we have also used the equation
\begin{equation*}
\int_C\Phi\frac{dG}{G}=0,
\end{equation*}
which can be easily proved as follows. Let $\Psi$ be a holomorphic function on
$\PP^1\setminus\Omega$ such that
\begin{equation*}
\Psi^\prime(z)=\frac{\partial\Phi(z)}{\partial z}~\text{and}~\Psi(\infty)=0.
\end{equation*}
Such function exists since $\partial\Phi(z)/\partial z=O(z^{-2})$ as
$z\rightarrow\infty$. We have $\Phi=\Psi+\overline{\Psi}$ on $\PP^1\setminus\Omega$
and by Cauchy theorem
\begin{equation*}
\int_C\Psi\frac{dG}{G}=0.
\end{equation*}
Using the same argument as for the case $n>0$ above shows that the integral with
$\overline{\Psi}$ also vanishes.

We summarize these results as the following statement (cf.~\cite{K-K-MW-W-Z}).

\begin{theorem} The normalized $1$-point current correlation functions of free bosons on
$\PP^1$ parameterized by $C\in\Tilde{\mathcal{C}}_a$ for every $a>0$ satisfy the Ward
identities, given by the following Laurent expansions at $z=\infty$
\begin{equation*}
\langle\langle\jmath(z)\rangle\rangle=-\sum_{n=1}^{\infty}z^{-n-1}\frac{\partial\log\tau}
{\partial t_n}dz=\mathbf{d}^\prime \log\tau,
\end{equation*}
and
\begin{equation*}
\langle\langle\bar\jmath(z)\rangle\rangle=-\sum_{n=1}^{\infty}\bar{z}^{-n-1}\frac{\partial
\log\tau}{\partial \bar{t}_n}d\bar{z} =\mathbf{d}^{\prime\prime}\log\tau.
\end{equation*}
\end{theorem}
\begin{remark} Since $\langle\jmath(z)\rangle_{DBC}=\langle\bar\jmath(z)\rangle_{DBC}=0$ the
theorem can be also trivially interpreted as computing in terms of the $\tau$-function
the difference between 1-point correlation functions of free bosons on $\PP^1$
parameterized by $C$ and of free bosons on $\PP^1\setminus\Omega$ with the Dirichlet
boundary condition. This will be relevant for the 2-point correlation functions.
\end{remark}
\begin{corollary} (\cite{MWZ,W-Z,K,K-K-MW-W-Z}) The function $\log\tau\in
C^\infty(\mathcal{C},\RR)$ is a generating function for the harmonic moments of
interior:
\begin{equation*}
v_0=\frac{\partial\log\tau}{\partial t_0}~\text{and}~
v_n=\frac{\partial\log\tau}{\partial t_n},~n\in\NN.
\end{equation*}
\end{corollary}
\begin{remark} It is instructive to compare this corollary with our results
with P.G.~Zograf~\cite{Z-T1,Z-T2} on the solution of the Poincar\'{e} problem of
accessory parameters (see also~\cite{Tak2} for an overview). Namely, let $X:=
\PP^1\setminus\{z_1,\dotsc,z_n\}$ be an $n$-punctured Riemann sphere, $n>3$,
normalized by $z_{n-2}=0, z_{n-1}=1~\text{and}~z_n=\infty$, and let $J: \HH\rightarrow
X$ be the uniformization map --- a complex-analytic covering of $X$ by the upper
half-plane $\HH$. The Schwarzian derivative $\mathcal{S}(J^{-1})$ of the inverse map
$J^{-1}$ is the following rational function on $\PP^1$
\begin{equation*}
\mathcal{S}(J^{-1})(z)=\sum_{i=1}^{n-3}\bigg(\frac{1/2}{(z-z_i)^2}+\frac{c_i}{z-z_i}\bigg)
+\frac{2-n}{2z(z-1)}.
\end{equation*}
The coefficients $c_i,\,i=1,\dotsc,n-3$, are smooth functions on the space of
punctures
\begin{equation*}
\mathcal{Z}_n:=\{(z_1,\dotsc, z_{n-3})\in\CC^{n-3}\,|\,z_i\neq z_j~\text{for}~i\neq j
~\text{and}~z_i\neq 0,1\}
\end{equation*}
and are called accessory parameters of the Fuchsian uniformization of
$n$-punctured Riemann spheres. We proved in~\cite{Z-T1,Z-T2} that, in accordance with a
conjecture of A.~Polyakov~\cite{P3}, there exists a real-valued smooth function $S$ on
$\mathcal{Z}_n$ such that
\begin{equation*}
c_i=-\frac{1}{2\pi}\frac{\partial S}{\partial z_i},~i=1,\dotsc, n-3.
\end{equation*}
The function $S$ is the critical value of the action function for the \emph{Liouville
theory --- the two-dimensional quantum gravity}, and the formulas for accessory
parameters follow from the \emph{semi-classical Ward identity for the $1$-point
correlation function with holomorphic component of the stress-energy tensor}
(see~\cite{Tak} for the details).

This comparison shows similarity between quantum theory of free bosons on $\PP^1$ in
the presence of the contour $C$ and quantum Liouville theory. Namely, 1-point Ward
identities for both theories imply that the logarithm of the $\tau$-function and the
critical value of the Liouville action (the logarithm of the semi-classical
approximation to the partition function) are generating functions for the harmonic
moments of interior and accessory parameters correspondingly.
\end{remark}
\begin{remark}
In the case when the punctures $z_1,\dotsc, z_{n-3}$ are real, the Riemann surface $X$
from the previous remark possesses an anti-holomorphic involution $z\mapsto \bar{z}$.
It is a classical result (see, e.g.~\cite{C-H}) that for this case the map $J^{-1}$ is
given by the conformal mapping of the upper half-plane $\HH$ (or the unit disk $D$)
onto a circular $n$-sided polygon with zero angles, inscribed into $S^1$. Under this
map the marked points $z_1,\dotsc, z_{n-3}, 0, 1,\infty$ on the boundary
$\RR\cup\{\infty\}$ of $\HH$ (or corresponding $n$ marked points on $S^1$) are mapped
onto the vertices of the polygon and this map is unique if the last three vertices are
normalized as $-1,-i,1$. The harmonic moments of the boundary of the polygon depend on
$n-3$ real parameters $z_1,\dotsc, z_{n-3}$. It would be instructive to express
accessory parameters through harmonic moments and compare results~\cite{Z-T1,Z-T2}
with the Corollary 3.3 directly. In order to get the analog of the Corollary 3.3 for
this case, one needs a generalization to the case of piece-wise smooth contours like
boundaries of circular polygons. We do not address this interesting question here.
\end{remark}
\begin{remark} As it was pointed out in the Remark 1.1 there is an inclusion
of $\Diff_+(S^1)/\Mob(S^1)$ into the Bers universal Teichm\"{u}ller space $T(1)$,
which contains all Teichm\"{u}ller spaces $T_{g,n}$ of Riemann surfaces of type
$(g,n)$ as complex submanifolds. Similar to the smooth case, there is a correspondence
$\gamma\mapsto C$, where $\gamma$ is a quasi-conformal homeomorphism of $S^1$ and $C$
is a quasi-circle --- an image of $S^1$ under quasi-conformal homeomorphism of the
complex plane $\CC$ which is conformal outside $S^1$. Extension of the above
formulation from smooth contours to quasi-circles would naturally allow to consider
Riemann surfaces of type $(g,n)$ by the same method. We do not address this important
question here.
\end{remark}

For completeness, let us show how to determine conformal map $G$ ``explicitly''~
\cite{MWZ,W-Z,K-K-MW-W-Z} from the Ward identity proved above. It follows from the
proof of the Krichever's lemma and the equation
\begin{equation*}
S_-(z)=\frac{\partial\Phi(z)}{\partial z} - \frac{t_0}{z}
\end{equation*}
that
\begin{equation*}
\frac{G^\prime(z)}{G(z)}=-\frac{\partial S_-(z)}{\partial t_0}= -\frac{\partial^2
\Phi(z)}{\partial t_0\partial z} + \frac{1}{z}.
\end{equation*}
Therefore
\begin{equation*}
\frac{G^\prime(z)}{G(z)}=\frac{1}{z} +
\sum_{n=1}^{\infty}z^{-n-1}\frac{\partial^2\log\tau} {\partial t_0 \partial t_n},
\end{equation*}
and integrating
\begin{equation*}
\log G(z) =\log b_{-1} + \log z -
\sum_{n=1}^{\infty}\frac{z^{-n}}{n}\frac{\partial^2\log\tau} {\partial t_0\partial
t_n}.
\end{equation*}
It is also possible~\cite{W-Z} to express the Robin's constant $\log b_{-1}=V$ (see
Section 1.2) through $\log\tau$. Namely, consider
\begin{equation*}
\frac{\partial^2\log\tau}{\partial t_0^2}=\frac{\partial v_0}{\partial t_0}.
\end{equation*}
Using the definition of $v_0$ and the calculus formula, we get
\begin{equation*}
\frac{\partial v_0}{\partial t_0}=\frac{1}{2\pi i}\int_C\log|z|^2\frac{dG}{G}=-2\log
b_{-1},
\end{equation*}
as can be readily shown by integration by parts. Therefore, one gets the result
in~\cite{W-Z,K-K-MW-W-Z}.
\begin{corollary}
The conformal map $G$ is given by the ``explicit formula''
\begin{equation*}
\log G(z) =\log z - \frac{1}{2}\frac{\partial^2\log\tau}{\partial t_0^2} -
\sum_{n=1}^{\infty}\frac{z^{-n}}{n}\frac{\partial^2\log\tau}{\partial t_0\partial
t_n}.
\end{equation*}
\end{corollary}
\begin{remark}
In complex analysis there is the following relation between moments $M_n$ of the
equilibrium distribution for the domain $\Omega$ and the conformal mapping $G$ of the
exterior domain $\PP^1\setminus\Omega$ (see, e.g.~\cite{Hille})
\begin{equation*}
\frac{G^\prime(z)}{G(z)}=\sum_{n=0}^\infty M_n z^{-n-1}.
\end{equation*}
As it follows for the above formula for $G^\prime(z)/G(z)$,
\begin{equation*}
M_n=\frac{\partial^2\log\tau}{\partial t_n\partial t_0},\,n\in\NN,
\end{equation*}
so that the smooth function $v_0=\partial\log\tau/\partial t_0$ on $\mathcal{C}$ is a generating
function for the moments $M_n$.
\end{remark}

Next, consider the normalized reduced 2-point current correlation function
\begin{equation*}
\langle\langle\jmath(z)\jmath(w)\rangle\rangle:=
\frac{\langle\jmath(z)\jmath(w)\rangle}{\langle\mathbf{1}\rangle_C}
-\langle\langle\jmath(z)\rangle\rangle\langle\langle\jmath(w)\rangle\rangle.
\end{equation*}
Using the same arguments as in the proof of Proposition 2.1: the change of variables
$\varphi=\Phi^{(\epsilon)}+ \tilde\varphi$ and passage to the limit
$\epsilon\rightarrow 0$, as well as the standard fact that
$\langle\tilde{\jmath}(z)\rangle_0=0$ we obtain, as in Section 4.1 of Part 1,
\begin{equation*}
\langle\jmath(z)\jmath(w)\rangle=\langle\tilde{\jmath}(z)\tilde{\jmath}(w) \rangle_0 +
\langle\mathbf{1}\rangle_C\frac{\partial \Phi(z)}{\partial z} \frac{\partial
\Phi(w)}{\partial w}dz\otimes dw.
\end{equation*}
As the result
\begin{align*}
\langle\langle\jmath(z)\jmath(w)\rangle\rangle & =\langle\langle\tilde{\jmath}(z)
\tilde{\jmath}(w)\rangle\rangle_0 =\pi \frac{\partial^2 G(z,w)}{\partial z\partial w}
dz\otimes dw \\ &= -\frac{dz\otimes dw}{(z-w)^2},
\end{align*}
where we used a simple expression for the Schiffer kernel $S$ on $\PP^1$ from Section
1.3.

Our goal is to compare this correlation function with the corresponding reduced
normalized 2-point correlation function for free bosons on $\PP^1\setminus\Omega$ with
the Dirichlet boundary condition. Since $\langle\jmath(z)\rangle_{DBC}=0$ we have
\begin{equation*}
\langle\langle\jmath(z)\jmath(w)\rangle\rangle_{DBC}:=
\frac{\langle\jmath(z)\jmath(w)\rangle_{DBC}}{\langle\mathbf{1}\rangle_{DBC}}.
\end{equation*}
Arguing as in the Section 4.1 of Part 1 and using results of Section 1.3 we get
\begin{equation*}
\langle\langle\jmath(z)\jmath(w)\rangle\rangle_{DBC}=\pi \frac{\partial^2
G_{DBC}(z,w)} {\partial z\partial w}dz\otimes dw=
-\frac{G^\prime(z)G^\prime(w)}{(G(z)-G(w))^2} dz\otimes dw.
\end{equation*}

On the other hand, consider $\partial^2\log\tau/\partial t_m\partial t_n,~m,n>0$.
Using the Ward identity for the 1-point correlation function, holomorphy of
$\partial\Phi/\partial z$ on $\PP^1\setminus\Omega$ and the calculus formula, we get
\begin{equation*}
\frac{\partial^2\log\tau}{\partial t_m\partial t_n}=-\frac{1}{2\pi
i}\frac{\partial}{\partial t_m} \int_C\frac{\partial \Phi(z)}{\partial z} z^n dz =
-\frac{1}{2\pi i}\int_C \frac{\partial^2\Phi(z)}{\partial t_m\partial z} z^n dz.
\end{equation*}
Using $\partial\Phi(z)/\partial z=S_{-}(z)+t_0/z$ and the formula
\begin{equation*}
\dot{\omega}^{(m)}_C=dz^m - \frac{\partial S_{-}}{\partial
t_m}dz,\,z\in\PP^1\setminus\Omega
\end{equation*}
(see the proof of part (i) of the Krichever's lemma), we have
\begin{equation*}
\frac{\partial^2\log\tau}{\partial t_m\partial t_n} =
 \frac{1}{2\pi i}\int_C w^n (\dot{\omega}^{(m)}_C-dw^m)=\frac{1}{2\pi i}
 \int_C w^n d(F_m \circ G - w^m),
\end{equation*}
where in the last equation we have used Krichever's lemma again. Next, it follows from
the definition of Faber polynomials (see Section 1.3) that
\begin{equation*}
\frac{G^\prime(z)G^\prime(w)}{(G(z)-G(w))^2}=\sum_{m=0}^\infty \frac{d F_m(G(w))}{dw}
z^{-m-1},
\end{equation*}
where $|G(z)|>|G(w)|$. From here we get
\begin{equation*}
\frac{d F_m(G(w))}{dw}=\frac{1}{2\pi
i}\int_{C_w}\frac{G^\prime(z)G^\prime(w)}{(G(z)-G(w))^2}z^m dz,
\end{equation*}
where contour $C_w$ is such that its image under the map $G$ contains the circle of
radius $|G(w)|$ inside. Similarly, from the expansion
\begin{equation*}
\frac{1}{(z-w)^2}=\sum_{m=0}m w^{m-1} z^{-m-1},
\end{equation*}
where $|z|>|w|$, we get
\begin{equation*}
m w^{m-1}=\frac{1}{2\pi i}\int_{|z|=R}\frac{1}{(z-w)^2}z^m dz,
\end{equation*}
where $|R|>|w|$. Thus we obtain
\begin{equation*}
\frac{d F_m(G(w))}{dw}-m w^{m-1}=\frac{1}{2\pi
i}\int_{C}\bigg(\frac{G^\prime(z)G^\prime(w)}
{(G(z)-G(w))^2}-\frac{1}{(z-w)^2}\bigg)z^m dz,
\end{equation*}
where we moved the contour of integration from the neighborhood of $\infty$ to $C$
because the integrand is regular for all $z,w\in\CC\setminus\Omega$. Thus we finally
get
\begin{align*}
\frac{\partial^2\log\tau}{\partial t_m\partial t_n}&=\frac{1}{(2\pi
i)^2}\int_C\int_C\bigg(
\frac{G^\prime(z)G^\prime(w)}{(G(z)-G(w))^2}-\frac{1}{(z-w)^2}\bigg)z^m w^n dz dw\\
&=\frac{1}{(2\pi i)^2}\int_C\int_Cz^m
w^n\left(\langle\langle\jmath(z)\jmath(w)\rangle\rangle -
\langle\langle\jmath(z)\jmath(w)\rangle\rangle_{DBC} \right)
\end{align*}
--- the Ward identity for the difference between normalized reduced 2-point correlation
functions of holomorphic current components for free bosons on $\PP^1$ parameterized
by $C$ and free bosons on $\PP^1\setminus\Omega$ with the Dirichlet boundary
condition.

Next, consider normalized reduced 2-point correlation function of holomorphic and
anti-holomorphic current components
\begin{equation*}
\langle\langle\jmath(z)\bar{\jmath}(w)\rangle\rangle:=
\frac{\langle\jmath(z)\bar{\jmath}(w)\rangle}{\langle\mathbf{1}\rangle_C}
-\langle\langle\jmath(z)\rangle\rangle\langle\langle\bar{\jmath}(w)\rangle\rangle.
\end{equation*}
As before we get
\begin{equation*}
\langle\langle\jmath(z)\bar{\jmath}(w)\rangle\rangle=\langle\jmath(z)
\bar{\jmath}(w)\rangle_0=\pi \frac{\partial^2 G(z,w)}{\partial z\partial \bar{w}}=0,
\end{equation*}
since the Bergman reproducing kernel on $\PP^1$ is 0 (see Section 1.2). Similarly,
arguing as in Section 4.1 of Part 1 and using results in Section 1.2 on the Bergman
kernel on $\PP^1\setminus\Omega$, we get
\begin{equation*}
\langle\langle\jmath(z)\bar{\jmath}(w)\rangle\rangle_{DBC} =\pi \frac{\partial^2
G_{DBC}(z,w)}{\partial z\partial\bar{w}}dz\otimes d\bar{w} =
-\frac{G^\prime(z)\overline{G^\prime(w)}}{(1-G(z)\overline{G(w)})^2}dz\otimes
d\bar{w}.
\end{equation*}

Computation of $\partial^2\log\tau/\partial t_m\partial \bar{t}_n,~m,n>0$ is also
similar to the one done before. Namely, since the vector field
$\partial/\partial\bar{t}_n$ corresponds to the meromorphic $(1,0)$-form
$\dot{\omega}_C^{(\bar{n})}$ on $\PP^1_C$ which coincides with $ - d(\overline{F_n\circ
1/\overline{G}})$ on the domain $\PP^1\setminus\Omega$ (see the proof of Krichever's
lemma in Section 1.3), we get
\begin{align*}
\frac{\partial^2\log\tau}{\partial t_m\partial \bar{t}_n} & = -\frac{1}{2\pi
i}\frac{\partial}{\partial\bar{t}_n}\int_C z^m \frac{\partial\Phi(z)}{\partial z}dz =
-\frac{1}{2\pi i}\int_C z^m \frac{\partial S_-(z)}{\partial\bar{t}_n}dz \\ & =
\frac{1}{2\pi i}\int_C z^m \dot{\omega}_C^{(\bar{n})} = -\frac{1}{2\pi i}\int_C z^m
d(\overline{F_n \circ 1/\overline{G}}).
\end{align*}
Using the representation
\begin{equation*}
-\frac{G^\prime(z)\overline{G^\prime(w)}}{(1-G(z)\overline{G(w)})^2}=\sum_{n=0}^\infty
\frac{d\overline{F_n(1/\overline{G(z)})}}{dz}\bar{w}^{-n-1},
\end{equation*}
which follows from the definition of Faber polynomials, we get by the Cauchy integral
formula
\begin{equation*}
\frac{d\overline{F_n(1/\overline{G(z)})}}{dz}=\frac{1}{2\pi i}\int_{C_+}
\frac{G^\prime(z)\overline{G^\prime(w)}}{(1-G(z)\overline{G(w)})^2}\bar{w}^nd\bar{w},
\end{equation*}
where $C_+$ is an arbitrary contour around 0 that contains the contour $C$ inside. As
a result we finally obtain
\begin{align*}
\frac{\partial^2\log\tau}{\partial t_m\partial \bar{t}_n} & = -\frac{1}{(2\pi
i)^2}\int_{C_+}\int_{C_+}
\frac{G^\prime(z)\overline{G^\prime(w)}}{(1-G(z)\overline{G(w)})^2}z^m \bar{w}^n dz
d\bar{w} \\ & = -\frac{1}{(2\pi i)^2}\int_{C_+}\int_{C_+} z^m \bar{w}^n
\left(\langle\langle\jmath(z)\bar{\jmath}(w)\rangle\rangle -
\langle\langle\jmath(z)\bar{\jmath}(w)\rangle\rangle_{DBC} \right)
\end{align*}
--- the Ward identity for the 2-point correlation function with holomorphic and
anti-holomorphic current components that computes the difference between correlation
functions of free bosons on $\PP^1$ parameterized by $C$ and of free bosons on
$\PP^1\setminus\Omega$ with the Dirichlet boundary condition. We summarize these
results in the following statement (cf.~\cite{K-K-MW-W-Z}).

\begin{theorem} Normalized reduced $2$-point current correlation functions for free bosons
on $\PP^1$ parameterized by $C\in\Tilde{\mathcal{C}}_a$ for every $a>0$ satisfy the
Ward identities, given by the following Laurent series expansions at $z=w=\infty$
\begin{align*}
\langle\langle\jmath(z)\jmath(w)\rangle\rangle -
\langle\langle\jmath(z)\jmath(w)\rangle\rangle_{DBC}&=
\left(\frac{G^\prime(z)G^\prime(w)}{(G(z)-G(w))^2} - \frac{1}
{(z-w)^2}\right)dz\otimes dw \\ &=\sum_{m,n=1}^{\infty} z^{-m-1}
w^{-n-1}\frac{\partial^2\log\tau}{\partial t_m\partial t_n} dz\otimes dw
\end{align*}
and
\begin{align*}
\langle\langle\jmath(z)\bar{\jmath}(w)\rangle\rangle -
\langle\langle\jmath(z)\bar{\jmath}(w)\rangle\rangle_{DBC} &=
\frac{G^\prime(z)\overline{G^\prime(w)}}{(1-G(z)\overline{G(w)})^2}dz\otimes d\bar{w}
\\ &=\sum_{m,n=1}^{\infty}z^{-m-1}\bar{w}^{-n-1}\frac{\partial^2\log\tau}{\partial
t_m\partial\bar{t}_n} dz\otimes d\bar{w} \\
&=\mathbf{d}^\prime\mathbf{d}^{\prime\prime}\log\tau.
\end{align*}
All higher reduced multi-point current correlation functions vanish.
\end{theorem}
\begin{corollary} For every $a>0$ he Hermitian metric $H$ on $\Tilde{\mathcal{C}}_a$ is
K\"{a}hler with the K\"{a}hler potential $\log\tau$.
\end{corollary}
\begin{proof} Immediately follows from the definition of the metric $H$ in Section 1.3.
\end{proof}
\begin{remark} This corollary should be compared
with the result~\cite{Z-T1,Z-T2,Z-T3} that the critical value $S$ of the Liouville
action functional is a potential for the Weil-Petersson K\"{a}hler metric on the space
of punctures $\mathcal{Z}_n$ (and, therefore, on the Teichm\"{u}ller space $T_{0,n}$)
\begin{equation*}
g^{i\bar{j}}_{WP}=-\frac{\partial^2 S}{\partial
z_i\partial\bar{z}_j},\,i,j=1,\dotsc,n-3.
\end{equation*}
As was noted in the Remark 1.7, the metric $H$ on $\Tilde{\mathcal{C}}_a$ is a
simplified analog of the Weil-Petersson metric on Teichm\"{u}ller spaces $T_{g,n}$.
Theorem 3.9 expresses it in terms of the 2-point current correlation function with
holomorphic and anti-holomorphic components. Similarly, the Weil-Petersson metric on
$T_{g,n}$ can be characterized as a semi-classical limit of a 2-point correlation
function with holomorphic and anti-holomorphic components of the stress-energy
tensor~\cite{Tak}.
\end{remark}

As before, it is possible to integrate explicitly the formulas in the theorem 3.9,
using the formula for $\log b_{-1}$. We obtain the following
result~\cite{W-Z,K-K-MW-W-Z}.
\begin{corollary}
\begin{align*}
&\log\frac{G(z)-G(w)}{z-w}=-\frac{1}{2}\frac{\partial^2\log\tau}{\partial t_0^2}
+\sum_{m,n=1}^{\infty}\frac{z^{-m}w^{-n}}{mn}\frac{\partial^2\log\tau}{\partial
t_m\partial t_n} \\ \intertext{and} \\ &
\log\left(\frac{G(z)\overline{G(w)}}{G(z)\overline{G(w)}-1}\right)=
\sum_{m,n=1}^{\infty}\frac{z^{-m}\bar{w}^{-n}}{mn}\frac{\partial^2\log\tau}{\partial
t_m\partial\bar{t}_n}.
\end{align*}
\end{corollary}

\end{document}